\documentclass[11pt]{article}

\usepackage{amsmath}
\usepackage{amssymb}
\usepackage{amsthm}
\usepackage{indentfirst}
\usepackage{enumerate}
\usepackage{cite}
\usepackage{bigints}

\newtheorem{thm}{Theorem}[section]

\newtheorem{cor}{Corollary}[section]
\newtheorem{rem}{Remark}[section]

\numberwithin{equation}{section}

\title{\bf LINEAR INDEPENDENCE\\
of Covariant Derivatives and Space-Curvatures}
\author{Nenad O. Vesi\'c\footnote{Faculty of Science and Mathematics, Ni\v s, Serbia,
Serbian Ministry of Education, Science and Technological
Development, Grant No. 174012}}
\date{}

\makeatletter
\def\maketag@@@#1{\hbox{\m@th\normalfont\normalsize#1}}
\makeatother

\begin{document}
  \maketitle

  \begin{abstract}
    The considerations about curvature tensors and pseudotensors for a non-symmetric
    affine connection space \big(see S. M. Min\v ci\'c, \cite{mincic2,
    mincic4,mincicnovi,mincvel}, M. Prvanovi\'c \cite{mileva1}\big) are advanced in this paper.
    It is obtained which kinds of  covariant
    derivatives, and identities of Ricci Type and curvatures
    of a non-symmetric affine connection space as well, are linearly independent.
     At the end of this research, curvature tensors of a special generalized
    Riemannian space are physically interpreted.
    \\[5pt]

    \textbf{Key words:} linear independence, affine connection, curvature tensor,
    energy-momentum tensor\\[2pt]

    \textbf{Math. Subj. Clasiffication:} 53A55,
    53B05, 53C15, 53Z05
  \end{abstract}

  \section{Introduction}

  An $N$-dimensional manifold $\mathcal M_N$ equipped with an affine
  connection $\nabla$ (with torsion) is the generalized affine connection
  space $\mathbb{GA}_N$ \big(see \cite{mincic2, mincic4, eisNRG, mincvel, mincicnovi, z4, mileva1}\big).
  The affine connection coefficients of the
  affine connection $\nabla$ are $L^i_{jk}$ and they are
  non-symmetric by indices $j$ and $k$.

  The symmetric and anti-symmetric parts of the coefficients
  $L^i_{jk}$ are

    \begin{eqnarray}
    L^i_{\underline{jk}}=\frac12\big(L^i_{jk}+L^i_{kj}\big)&\mbox{and}&
    L^i_{\underset\vee{jk}}=\frac12\big(L^i_{jk}-L^i_{kj}\big).
    \label{eq:Lsimantisim}
  \end{eqnarray}

  The symmetric parts $L^i_{\underline{jk}}$ are the affine connection
  coefficients of the corresponding symmetric affine connection $\overset0\nabla$
  because they satisfy the corresponding transformation
  rule \cite{mik5,sinjukov}. The affine connection
  space $\mathbb A_N$  equipped with affine connection $\overset0\nabla$
   is the associated space (of the
  space $\mathbb{GA}_N$). The anti-symmetric part
  $L^i_{\underset\vee{jk}}$ is equal to a half of the torsion tensor
  for the space $\mathbb{GA}_N$.

  \subsection{Symmetric affine connection space}

  As we said above, the $N$-dimensional manifold $\mathcal M_N$
  equipped with a symmetric affine connection $\overset0\nabla$
   is the symmetric affine connection space $\mathbb A_N$
  (see \cite{mik5,sinjukov}). One kind of covariant
  differentiation with respect to the affine connection
  $\overset0\nabla$ exists

  \begin{equation}
    a^i_{j|k}=a^i_{j,k}+L^i_{\underline{\alpha k}}a^\alpha_j-
    L^\alpha_{\underline{jk}}a^i_\alpha,
    \label{eq:covderivativesim}
  \end{equation}

  \noindent In the last equation, the object $a^i_j$ is a tensor of
  the type $(1,1)$ and the partial derivative $\partial/\partial
  x^k$ is denoted by comma.

  It is founded one Ricci Type identity with respect to the affine
  connection $\overset0\nabla$. Moreover, it exists
  one curvature tensor of the space $\mathbb A_N$

  \begin{equation}
    R^i_{jmn}=L^i_{\underline{jm},n}-L^i_{\underline{jn},m}+
    L^\alpha_{\underline{jm}}L^i_{\underline{\alpha n}}-
    L^\alpha_{\underline{jn}}L^i_{\underline{\alpha m}},
    \label{eq:R}
  \end{equation}

  Many authors have promoted the theory of symmetric affine connection spaces.
  Some of them are J. Mike\v s and his research group \cite{mik5},
  N. S. Sinyukov \cite{sinjukov} and many others.

  \subsection{Non-symmetric affine connection space}

  An $N$-dimensional manifold equipped with a non-symmetric affine connection
  $\nabla$ is the non-symmetric affine connection space
  $\mathbb{GA}_N$ \big(see \cite{cvetkoviczlatanovic1, z4, mincic2, mincic4, mincicnovi, mincvel,
  mileva1}\big).
  Unlike in the case of
  a symmetric affine connection space, there are defined four kinds of
  covariant derivatives \cite{mincic2, mincic4, mileva1, mincicnovi, mincvel, z4,
  cvetkoviczlatanovic1}

    \begin{equation}
    \begin{array}{ll}
      a^i_{j\underset1|k}=a^i_{j,k}+L^i_{\alpha
      k}a^\alpha_j-L^\alpha_{jk}a^i_\alpha,&
      a^i_{j\underset2|k}=a^i_{j,k}+L^i_{k\alpha}a^\alpha_j
      -L^\alpha_{kj}a^i_\alpha,\\
      a^i_{j\underset3|k}=a^i_{j,k}+L^i_{\alpha
      k}a^\alpha_j-L^\alpha_{kj}a^i_\alpha,&
      a^i_{j\underset4|k}=a^i_{j,k}+L^i_{k\alpha}a^\alpha_j
      -L^\alpha_{jk}a^i_\alpha.
    \end{array}
    \label{eq:covderivativensim}
  \end{equation}

  With respect to these kinds of covariant derivatives, it is obtained
  four curvature tensors, eight derived curvature tensors and
  fifteen curvature pseudotensors \cite{mincic2, mincic4, mincicnovi, mincvel,
  cvetkoviczlatanovic1,z4} of the
  space $\mathbb{GA}_N$. In the papers about non-symmetric
  affine connection spaces, authors have cited and
  confirmed the following statements \cite{mincic2, mincic4, mincvel}

  \begin{enumerate}[Statement 1:]
    \item[Statement 1:] \emph{In the set of twelve curvature tensors of the space
    $\mathbb{GA}_N$ with non-symmetric affine connection $\nabla$,
    there are five independent ones, while the rest can be expressed
    in terms of these five tensors and the curvature tensor of the associated space.}
    \item[Statement 2:] \emph{All fifteen curvature pseudotensors are given in the
    terms of five of them, five curvature tensors of the space
    $\mathbb{GA}_N$ and the curvature tensor of the associated space
    $\mathbb{A}_N$.}
    \item[Statement 3:] \emph{Curvature pseudotensors have the role
    and the form of the curvature tensors}.
  \end{enumerate}

  The role of the curvature tensor (\ref{eq:R}) of the symmetric
  affine connection space $\mathbb{A}_N$ is to \emph{uniquely}
  generate the curvature $a^i_{j|m|n}-a^i_{j|n|m}$ of this space.
  In other words, the geometrical object $a^\alpha_jR^i_{\alpha
  mn}-a^i_\alpha R^\alpha_{jmn}$ is the only component of this
  curvature.

  In this paper, we are aimed to find all components for the
  curvatures $a^i_{j\underset p|m\underset q|n}-
  a^i_{j\underset r|n\underset s|m}$, $p,q,r,s\in\{1,2,3,4\}$, of the space
  $\mathbb{GA}_N$. Moreover, we will obtain and prove what are roles
  of the curvature tensors, torsion tensor and curvature
  pseudotensors with respect to the curvatures of the space
  $\mathbb{GA}_N$ in this study.

  \subsubsection{Generalized Riemannian space}

  Special kind of affine connection spaces with torsion are the
  generalized Riemannian spaces $\mathbb{GR}_N$ \big(in the sense of
  Eisenhart \cite{eis01,eis02}\big). These spaces are manifolds $\mathcal M_N$
  equipped with a non-symmetric metric tensor $g_{ij}$.

  The symmetric and anti-symmetric part of the tensor $g_{ij}$ are

  \begin{eqnarray}
    g_{\underline{ij}}=\frac12(g_{ij}+g_{ji})&\mbox{and}&
    g_{\underset\vee{ij}}=\frac12(g_{ij}-g_{ji}).
    \label{eq:gsimantisim}
  \end{eqnarray}

  We assume that the matrix $\big[g_{\underline{ij}}\big]_{N\times
  N}$ is non-singular. The contravariant symmetric metric tensor
  $g^{\underline{ij}}$ is obtained by the inverting the matrix
  $\big[g_{\underline{ij}}\big]$. Hence, it holds the equality
  $g^{\underline{i\alpha}}g_{\underline{j\alpha}}=\delta^i_j$, for
  the Kronecker delta-symbols $\delta^i_j$.

  The affine connection coefficients of the space $\mathbb{GR}_N$
  are the generalized Christoffel symbols \cite{eis01, eis02}

  \begin{equation}
    \Gamma^i_{jk}=\frac12g^{\underline{i\alpha}}\big(g_{j\alpha,k}-g_{jk,\alpha}+
    g_{\alpha k,j}\big).
    \label{eq:christoffelsymbolgen}
  \end{equation}

  The symmetric part of the coefficient $\Gamma^i_{jk}$,
  $\Gamma^i_{\underline{jk}}=\dfrac12\big(\Gamma^i_{jk}+\Gamma^i_{kj}\big)$,
  is the corresponding Christoffel symbol. The anti-symmetric part
  of this symbol equal to a half of the torsion tensor of the
  space $\mathbb{GR}_N$, i.e.
  $\Gamma^i_{\underset\vee{jk}}=\dfrac12\big(\Gamma^i_{jk}-\Gamma^i_{kj}\big)$,
  is

  \begin{equation}
    \Gamma^i_{\underset\vee{jk}}
    =\frac12g^{\underline{i\alpha}}\big(g_{\underset\vee{j\alpha},k}
    -g_{\underset\vee{jk},\alpha}+
    g_{\underset\vee{\alpha k},j}\big).
    \label{eq:christoffelsymbolgenantisim}
  \end{equation}

  \subsection{Motivation}

  A. Einstein \cite{e1,e2,e3} was the first who applied
  the complex metrics (equivalent with the non-symmetric metrics)
  in his research about the Unified Field Theory. The affine connection
  coefficients in Einstein's works are not explicit functions
  of non-symmetric metrics but they satisfy the Einstein Metricity Condition
  \begin{eqnarray}
    g_{ij,k}-\Gamma^\alpha_{ik}g_{\alpha
    j}-\Gamma^\alpha_{kj}g_{i\alpha}=0.
    \label{eq:einsteinmetricitycondition}
  \end{eqnarray}
  With
  respect to Einstein's work, the symmetric parts
  $\Gamma^i_{\underline{jk}}$ of the affine connection coefficients refer
  to gravity but the anti-symmetric parts $\Gamma^i_{\underset\vee{jk}}$
  are important for researching in the theory of electromagnetism.

   S. Ivanov and M. Lj.
  Zlatanovi\'c \cite{z2} contributed to the theory of
  differential geometry motivated by Einstein's work.

  Motivated by the Einstein's considerations \cite{e1, e2, e3} and
  Eisenhart's  results \cite{eisNRG, eis01, eis02},  many authors have
  developed the theory of
  non-symmetric affine connection spaces \cite{mileva1,
  mincic2,mincic4, mincicnovi, mincvel,cvetkoviczlatanovic1, z4}
  and many others. The identities of Ricci Type in their papers are
  confirmations of the corresponding results initially presented in \cite{mincic2, mincic4}.
  We are aimed to computationally complete the derivation of the Ricci Type
  identities in this paper.

  The main purposes of the following research are:

  \begin{enumerate}[a.]
    \item To examine how many covariant derivatives of the covariant
    derivatives
    (\ref{eq:covderivativesim}, \ref{eq:covderivativensim}) may be expressed
    as the linear combinations of the rest ones,
    \item To find linearly independent identities of Ricci type (linearly independent
    curvatures) for the space $\mathbb{GA}_N$,
    \item To find linearly independent curvatures of the space
    $\mathbb{GA}_N$,
     \item To interpret the curvature tensors of the space $\mathbb{GR}_4$
     with respect to cosmology.
  \end{enumerate}

  \section{Kinds of covariant derivatives: four plus one but three}

  By the equations (\ref{eq:covderivativesim},
  \ref{eq:covderivativensim}), the covariant derivatives are
  expressed in the corresponding vector forms. The summands at the
  right sides of the equalities in these equations are the
  components of the corresponding covariant derivatives (treated as
  vectors). Analogously as in the case of curvature tensors
  \cite{mincic2, mincic4, mincvel, mincicnovi, z4}, we are aimed to
  search how many of these vectors are linearly independent at the start of this
  section. After that, we will obtain the linearly independent
  identities of Ricci Type and the curvature tensors from these
  identities.

  With respect to the equations (\ref{eq:covderivativesim},
  \ref{eq:covderivativensim}), we get:

  \begin{align}
    &a^{i}_{j|k}=
    \frac12a^{i}_{j\underset1|k}+
    \frac12a^{i}_{j\underset2|k},\label{eq:0|=1|+2|}\\\displaybreak[0]
    &a^{i}_{j|k}=
    \frac12a^{i}_{j\underset3|k}+
    \frac12a^{i}_{j\underset4|k},\label{eq:0|=3|+4|}\\\displaybreak[0]
    &a^{i}_{j\underset1|k}=
    2a^{i}_{j|k}-a^{i}_{j\underset2|k},\label{eq:1|=0|+2|}\\\displaybreak[0]
    &a^{i}_{j\underset1|k}=
    -a^{i}_{j\underset2|k}+
    a^{i}_{j\underset3|k}+
    a^{i}_{j\underset4|k},\label{eq:1|=2|+3|+4|}\\\displaybreak[0]
    &a^{i}_{j\underset2|k}=
    2a^{i}_{j|k}-
    a^{i}_{j\underset1|k},\label{eq:2|=0|+1|}\\\displaybreak[0]
    &a^{i}_{j\underset2|k}=
    -a^{i}_{j\underset1|k}+
    a^{i}_{j\underset3|k}+
    a^{i}_{j\underset4|k},\label{eq:2|=1|+3|+4|}\\\displaybreak[0]
    &a^{i}_{j\underset 3|k}=2
    a^{i}_{j|k}-
    a^{i}_{j\underset4|k},\label{eq:3|=0|+4|}\\\displaybreak[0]
    &a^{i}_{j\underset3|k}=
    a^{i}_{j\underset1|k}+
    a^{i}_{j\underset2|k}-
    a^{i}_{j\underset4|k},\label{eq:3|=1|+2|+4|}\\\displaybreak[0]
    &a^{i}_{j\underset4|k}=
    2a^{i}_{j|k}-
    a^{i}_{j\underset3|k},\label{eq:4|=0|+3|}\\\displaybreak[0]
    &a^{i}_{j\underset4|k}=
    a^{i}_{j\underset1|k}+
    a^{i}_{j\underset2|k}-
    a^{i}_{j\underset3|k},\label{eq:4|=1|+2|+3|}
  \end{align}

  \noindent for the tensor $a^{i}_{j}$ of the type $(1,1)$.




  Based on the equations
  (\ref{eq:0|=1|+2|} - \ref{eq:4|=1|+2|+3|}), we conclude that the
  covariant derivatives $a^i_{j|k}$ and $a^i_{j\underset4|k}$ may be expressed as the linear
  combinations of the covariant derivatives $a^i_{j\underset1|k}$,
  $a^i_{j\underset2|k}$ and $a^i_{j\underset3|k}$. As the tensors,
  the covariant derivatives $a^i_{j\underset1|k}$,
  $a^i_{j\underset2|k}$ and $a^i_{j\underset3|k}$ are linearly
  independent what may be checked in the standard way.

  For this reason, three kinds of covariant derivatives whose
  corresponding linear combinations express the rest two kinds of
  covariant derivatives will be called \emph{the linearly
  independent} (kinds of) covariant derivatives.

  Therefore, it holds the following theorem:

  \begin{thm}
    Three of the covariant derivatives
    $a^i_{j|k}$, $a^i_{j\underset1|k}$, $a^i_{j\underset2|k}$,
    $a^i_{j\underset3|k}$, $a^i_{j\underset4|k}$
    are linearly independent.\hfill\qed
  \end{thm}

  \begin{cor}
The triples of linearly independent covariant derivatives
$a^i_{j|k}$, $a^i_{j\underset1|k}$, $a^i_{j\underset2|k}$,
    $a^i_{j\underset3|k}$, $a^i_{j\underset4|k}$ are

\begin{footnotesize}
    \begin{eqnarray*}
      \begin{array}{cccc}
        b_1=\Big\{a^i_{j\underset1|k},
        a^i_{j\underset2|k},a^i_{j\underset3|k}\Big\}
        ,&
        b_2=\Big\{a^i_{j\underset1|k},a^i_{j\underset2|k},a^i_{j\underset4|k}\Big\}
        ,&
        b_3=\Big\{a^i_{j\underset1|k},a^i_{j\underset3|k},a^i_{j\underset4|k}\Big\},&
        b_4=\Big\{a^i_{j\underset2|k},a^i_{j\underset3|k},a^i_{j\underset4|k}\Big\},\\
        b_5=\Big\{a^i_{j|k},a^i_{j\underset1|k},a^i_{j\underset3|k}\Big\},&
        b_6=\Big\{a^i_{j|k},a^i_{j\underset1|k},a^i_{j\underset4|k}\Big\},&
        b_7=\Big\{a^i_{j|k},a^i_{j\underset2|k},a^i_{j\underset3|k}\Big\},&
        b_8=\Big\{a^i_{j|k},a^i_{j\underset2|k},a^i_{j\underset4|k}\Big\},
      \end{array}
    \end{eqnarray*}
\end{footnotesize}

\noindent unlike the covariant derivatives in the triples
$\big\{a^i_{j|k},a^i_{j\underset1|k},a^i_{j\underset2|k}\big\}$ and
$\big\{a^i_{j|k},a^i_{j\underset3|k},a^i_{j\underset4|k}\big\}$.
\qed
  \end{cor}

  \subsection{Identities of Ricci Type}

  Ricci-Type identities with respect to non-symmetric affine
  connection are obtained in many papers. They are initially
  founded in \big(S. M. Min\v ci\'c, \cite{mincic2, mincic4}\big).
  After that, many authors confirmed these identities for different
  non-symmetric affine connection spaces \big(see the papers \cite{mincvel,mincicnovi,
  cvetkoviczlatanovic1,z4}\big).

  In the last cited papers, the authors combined the first and the
  second kind of covariant derivatives (\ref{eq:covderivativensim})
  together with one identity with respect to the third and the fourth kind
  of covariant derivatives.
  Because the covariant
  derivatives $\underset1|,\underset2|,\underset3|$ are linearly
  independent, we will complete the computations about
  Ricci-Type identities with respect to  the linearly independent
  covariant derivatives $\underset1|,\underset2|,\underset3|$.
  Moreover, we will simplify the obtained identities of
  Ricci Type obtained in
  \cite{mincic2, mincic4, mincicnovi, mincvel, z4,cvetkoviczlatanovic1}.

  The double
  covariant derivatives of the tensor
  $a^i_j$ of the type $(1,1)$ with respect to these three linearly
  independent ones are

  \begin{align}
    &\aligned
    a^i_{j\underset1|m\underset1|n}&=a^i_{j,mn}-L^\alpha_{jn}a^i_{\alpha,m}
    -L^\alpha_{jm}a^i_{\alpha,n}-L^\alpha_{mn}a^i_{j,\alpha}
    +L^i_{\alpha n}a^\alpha_{j,m}+
    L^i_{\alpha m}a^\alpha_{j,n}\\
    &+a^\alpha_j\big(
    L^i_{\alpha m,n}+L^\beta_{\alpha m}L^i_{\beta n}
    -L^i_{\alpha\beta}L^\beta_{mn}\big)
    -a^i_\alpha\big(
    L^\alpha_{jm,n}-L^\alpha_{\beta
    m}L^\beta_{jn}-L^\alpha_{j\beta}L^\beta_{mn}
    \big)\\&-
    a^\alpha_\beta\big(L^i_{\alpha m}L^\beta_{jn}+L^i_{\alpha n}L^\beta_{jm}\big),
    \endaligned\label{eq:||11}
    \end{align}
    \begin{align}
    &\aligned
    a^i_{j\underset1|m\underset2|n}&=a^i_{j,mn}-L^\alpha_{nj}a^i_{\alpha,m}
    -L^\alpha_{jm}a^i_{\alpha,n}-L^\alpha_{nm}a^i_{j,\alpha}
    +L^i_{n\alpha}a^\alpha_{j,m}+
    L^i_{\alpha m}a^\alpha_{j,n}\\
    &+a^\alpha_j\big(
    L^i_{\alpha m,n}+L^\beta_{\alpha m}L^i_{n\beta}
    -L^i_{\alpha\beta}L^\beta_{nm}\big)
    -a^i_\alpha\big(
    L^\alpha_{jm,n}-L^\alpha_{\beta
    m}L^\beta_{nj}-L^\alpha_{j\beta}L^\beta_{nm}
    \big)\\&-
    a^\alpha_\beta\big(L^i_{\alpha m}L^\beta_{nj}+L^i_{n\alpha}L^\beta_{jm}\big),
    \endaligned\label{eq:||12}
    \\\displaybreak[0]
    &\aligned
    a^i_{j\underset1|m\underset3|n}&=a^i_{j,mn}-L^\alpha_{nj}a^i_{\alpha,m}
    -L^\alpha_{jm}a^i_{\alpha,n}-L^\alpha_{nm}a^i_{j,\alpha}
    +L^i_{\alpha n}a^\alpha_{j,m}+
    L^i_{\alpha m}a^\alpha_{j,n}\\
    &+a^\alpha_j\big(
    L^i_{\alpha m,n}+L^\beta_{\alpha m}L^i_{\beta n}
    -L^i_{\alpha\beta}L^\beta_{nm}\big)
    -a^i_\alpha\big(
    L^\alpha_{jm,n}-L^\alpha_{\beta
    m}L^\beta_{nj}-L^\alpha_{j\beta}L^\beta_{nm}
    \big)\\&-
    a^\alpha_\beta\big(L^i_{\alpha m}L^\beta_{nj}+L^i_{\alpha n}L^\beta_{jm}\big),
    \endaligned\label{eq:||13}\\\displaybreak[0]
    &\aligned
    a^i_{j\underset2|m\underset1|n}&=a^i_{j,mn}-L^\alpha_{jn}a^i_{\alpha,m}
    -L^\alpha_{mj}a^i_{\alpha,n}-L^\alpha_{mn}a^i_{j,\alpha}
    +L^i_{\alpha n}a^\alpha_{j,m}+
    L^i_{m\alpha}a^\alpha_{j,n}\\
    &+a^\alpha_j\big(
    L^i_{m\alpha,n}+L^\beta_{m\alpha}L^i_{\beta n}
    -L^i_{\beta\alpha}L^\beta_{mn}\big)
    -a^i_\alpha\big(
    L^\alpha_{mj,n}-L^\alpha_{m\beta}L^\beta_{jn}
    -L^\alpha_{\beta j}L^\beta_{mn}
    \big)\\&-
    a^\alpha_\beta\big(L^i_{m\alpha}L^\beta_{jn}+L^i_{\alpha n}L^\beta_{mj}\big),
    \endaligned\label{eq:||21}\\\displaybreak[0]
    &\aligned
    a^i_{j\underset2|m\underset2|n}&=a^i_{j,mn}-L^\alpha_{nj}a^i_{\alpha,m}
    -L^\alpha_{mj}a^i_{\alpha,n}-L^\alpha_{nm}a^i_{j,\alpha}
    +L^i_{n\alpha}a^\alpha_{j,m}+
    L^i_{m\alpha}a^\alpha_{j,n}\\
    &+a^\alpha_j\big(
    L^i_{m\alpha,n}+L^\beta_{m\alpha}L^i_{n\beta}
    -L^i_{\beta\alpha}L^\beta_{nm}\big)
    -a^i_\alpha\big(
    L^\alpha_{mj,n}-L^\alpha_{m\beta}L^\beta_{nj}
    -L^\alpha_{\beta j}L^\beta_{nm}
    \big)\\&-
    a^\alpha_\beta\big(L^i_{m\alpha}L^\beta_{nj}+L^i_{n\alpha}L^\beta_{mj}\big),
    \endaligned\label{eq:||22}\\\displaybreak[0]
    &\aligned
    a^i_{j\underset2|m\underset3|n}&=a^i_{j,mn}-L^\alpha_{nj}a^i_{\alpha,m}
    -L^\alpha_{mj}a^i_{\alpha,n}-L^\alpha_{nm}a^i_{j,\alpha}
    +L^i_{\alpha n}a^\alpha_{j,m}+
    L^i_{m\alpha}a^\alpha_{j,n}\\
    &+a^\alpha_j\big(
    L^i_{m\alpha,n}+L^\beta_{m\alpha}L^i_{\beta n}
    -L^i_{\beta\alpha}L^\beta_{nm}\big)
    -a^i_\alpha\big(
    L^\alpha_{mj,n}-L^\alpha_{m\beta}L^\beta_{nj}-L^\alpha_{\beta j}L^\beta_{nm}
    \big)\\&-
    a^\alpha_\beta\big(L^i_{m\alpha}L^\beta_{nj}+L^i_{\alpha n}L^\beta_{mj}\big),
    \endaligned\label{eq:||23}\\\displaybreak[0]
    &\aligned
    a^i_{j\underset3|m\underset1|n}&=a^i_{j,mn}-L^\alpha_{jn}a^i_{\alpha,m}
    -L^\alpha_{mj}a^i_{\alpha,n}-L^\alpha_{mn}a^i_{j,\alpha}
    +L^i_{\alpha n}a^\alpha_{j,m}+
    L^i_{\alpha m}a^\alpha_{j,n}\\
    &+a^\alpha_j\big(
    L^i_{\alpha m,n}+L^\beta_{\alpha m}L^i_{\beta n}
    -L^i_{\alpha\beta}L^\beta_{mn}\big)
    -a^i_\alpha\big(
    L^\alpha_{mj,n}-L^\alpha_{m\beta}L^\beta_{jn}-L^\alpha_{\beta j}L^\beta_{mn}
    \big)\\&-
    a^\alpha_\beta\big(L^i_{\alpha m}L^\beta_{jn}+L^i_{\alpha n}L^\beta_{mj}\big),
    \endaligned\label{eq:||31}\\\displaybreak[0]
    &\aligned
    a^i_{j\underset3|m\underset2|n}&=a^i_{j,mn}-L^\alpha_{nj}a^i_{\alpha,m}
    -L^\alpha_{mj}a^i_{\alpha,n}-L^\alpha_{nm}a^i_{j,\alpha}
    +L^i_{n\alpha}a^\alpha_{j,m}+
    L^i_{\alpha m}a^\alpha_{j,n}\\
    &+a^\alpha_j\big(
    L^i_{\alpha m,n}+L^\beta_{\alpha m}L^i_{n\beta}
    -L^i_{\alpha\beta}L^\beta_{nm}\big)
    -a^i_\alpha\big(
    L^\alpha_{mj,n}-L^\alpha_{m\beta}L^\beta_{nj}-L^\alpha_{\beta j}L^\beta_{nm}
    \big)\\&-
    a^\alpha_\beta\big(L^i_{\alpha m}L^\beta_{nj}+L^i_{n\alpha}L^\beta_{mj}\big),
    \endaligned\label{eq:||32}\\\displaybreak[0]
    &\aligned
    a^i_{j\underset3|m\underset3|n}&=a^i_{j,mn}-L^\alpha_{nj}a^i_{\alpha,m}
    -L^\alpha_{mj}a^i_{\alpha,n}-L^\alpha_{nm}a^i_{j,\alpha}
    +L^i_{\alpha n}a^\alpha_{j,m}+
    L^i_{\alpha m}a^\alpha_{j,n}\\
    &+a^\alpha_j\big(
    L^i_{\alpha m,n}+L^\beta_{\alpha m}L^i_{\beta n}
    -L^i_{\alpha\beta}L^\beta_{nm}\big)
    -a^i_\alpha\big(
    L^\alpha_{mj,n}-L^\alpha_{m\beta}L^\beta_{nj}-L^\alpha_{\beta j}L^\beta_{nm}
    \big)\\&-
    a^\alpha_\beta\big(L^i_{\alpha m}L^\beta_{nj}+L^i_{\alpha
    n}L^\beta_{mj}\big).
    \endaligned\label{eq:||33}
  \end{align}

  We conclude that there are $3\cdot3\cdot 3\cdot 3=81$ identities of Ricci Type obtained from
  the differences
   $a^i_{j\underset
  p|m\underset q|n}-a^i_{j\underset r|n\underset s|m}$, for
  $p,q,r,s\in\{1,2,3\}$.

  After using the equality
  $L^i_{jk}=L^i_{\underline{jk}}+L^i_{\underset\vee{jk}}$, one gets

  \begin{eqnarray}
    L^i_{jm,n}=L^i_{\underline{jm},n}+L^i_{\underset\vee{jm},n}&\mbox{and}&
    L^i_{jk}L^a_{bc}=L^i_{\underline{jk}}L^a_{\underline{bc}}+
    L^i_{\underset\vee{jk}}L^a_{\underline{bc}}+
    L^i_{\underline{jk}}L^a_{\underset\vee{bc}}+
    L^i_{\underset\vee{jk}}L^a_{\underset\vee{bc}}.
    \label{eq:LsimantisimLLsimantisim}
  \end{eqnarray}

  We substituted the last equation into the all of 81 identities
  of Ricci Type and founded that the next theorem and its corollaries
  are satisfied. We will prove this theorem for the case of $p=r=1$ and $q=s=2$
  in the subsection \ref{subsectiontransformationctpct}. All other identities of the Ricci Type
  from the next theorem may be proved analogously.

  \begin{thm}
    The family of Ricci-Type identities with respect to the
        covariant derivatives $\underset1|$, $\underset2|$,
        $\underset3|$ is

        \begin{footnotesize}
        \begin{equation}
        \aligned
          a^i_{j\underset p|m\underset q|n}-
          a^i_{j\underset r|n\underset
          s|m}&=2\overset1c{}_1L^\alpha_{\underset\vee{jm}}a^i_{\alpha|n}
          +2\overset1c{}_2L^\alpha_{\underset\vee{jn}}a^i_{\alpha|m}
          +2\overset1c{}_3L^\alpha_{\underset\vee{mn}}a^i_{j|\alpha}
          +2\overset1c{}_4L^i_{\underset\vee{\alpha n}}a^\alpha_{j|m}
          +2\overset1c{}_5L^i_{\underset\vee{\alpha
          m}}a^\alpha_{j|n}\\&
          +a^\alpha_j\big(R^i_{\alpha mn}
          +\overset1c{}_6L^i_{\underset\vee{\alpha m}|n}
          +\overset1c{}_7L^i_{\underset\vee{\alpha n}|m}
          +\overset1c{}_8L^\beta_{\underset\vee{\alpha
          m}}L^i_{\underset\vee{\beta n}}
          +\overset1c{}_9L^\beta_{\underset\vee{\alpha
          n}}L^i_{\underset\vee{\beta m}}
          +2\overset1c{}_{10}L^i_{\underset\vee{\alpha
          \beta}}L^\beta_{\underset\vee{mn}}
          \big)\\&
          -a^i_\alpha
          \big(
          R^\alpha_{jmn}+\overset1c{}_{11}L^\alpha_{\underset\vee{jm}|n}
          +\overset1c{}_{12}L^\alpha_{\underset\vee{jn}|m}
          +\overset1c{}_{13}L^\alpha_{\underset\vee{\beta
          n}}L^\beta_{\underset\vee{jm}}
          +\overset1c{}_{14}L^\alpha_{\underset\vee{\beta
          m}}L^\beta_{\underset\vee{jn}}
          +2\overset1c{}_{15}L^\alpha_{\underset\vee{j\beta}}L^\beta_{\underset\vee{mn}}
          \big)\\&
          -2a^\alpha_\beta\big(
          \overset1c{}_{16}L^i_{\underset\vee{\alpha
          m}}L^\beta_{\underset\vee{jn}}
          +\overset1c{}_{17}L^i_{\underset\vee{\alpha
          n}}L^\beta_{\underset\vee{jm}}
          \big),
        \endaligned\label{eq:vricciti}
        \end{equation}
        \end{footnotesize}

\noindent        for $p,q,r,s\in\{1,2,3\}$, the covariant derivative
with respect to the symmetric affine connection denoted by $|$ and
the corresponding coefficients
        $\overset1c{}_k\in\{0,1,-1\}$, $k=1,\ldots,17$.\qed
  \end{thm}
    \begin{cor}
    Seventeen of the identities from the family
    \emph{(\ref{eq:vricciti})} are linearly independent. Some seventeen of
    these linearly independent Ricci-Type identities are given by
    the equations \emph{(\ref{eq:ric11-11}---\ref{eq:ric33-33})} in the Appendix I.
    All other Ricci-Type identities are the corresponding linear
    combinations of the linearly independent ones.\qed
  \end{cor}
  \begin{cor}
    The family \emph{(\ref{eq:vricciti})} and the families

    \begin{scriptsize}
        \begin{equation}
        \aligned
          a^i_{j\underset p|m\underset q|n}-
          a^i_{j\underset r|n\underset
          s|m}&=
          2\overset1c{}_1d^1_1L^\alpha_{\underset\vee{jm}}a^i_{\alpha\underset1|n}
        +2\overset1c{}_1d^2_1L^\alpha_{\underset\vee{jm}}a^i_{\alpha\underset2|n}
        +2\overset1c{}_1d^3_1L^\alpha_{\underset\vee{jm}}a^i_{\alpha\underset3|n}
          +2\overset1c{}_2d^1_2L^\alpha_{\underset\vee{jn}}a^i_{\alpha\underset1|m}
        +2\overset1c{}_2d^2_2L^\alpha_{\underset\vee{jn}}a^i_{\alpha\underset2|m}
        \\&
        +2\overset1c{}_2d^3_2L^\alpha_{\underset\vee{jn}}a^i_{\alpha\underset3|m}
        +2\overset1c{}_3d^1_3L^\alpha_{\underset\vee{mn}}a^i_{j\underset1|\alpha}
        +2\overset1c{}_3d^2_3L^\alpha_{\underset\vee{mn}}a^i_{j\underset2|\alpha}
        +2\overset1c{}_3d^3_3L^\alpha_{\underset\vee{mn}}a^i_{j\underset3|\alpha}
          +2\overset1c{}_4d^1_4L^i_{\underset\vee{\alpha n}}a^\alpha_{j\underset1|m}
          \\&
        +2\overset1c{}_4d^2_4L^i_{\underset\vee{\alpha n}}a^\alpha_{j\underset2|m}
        +2\overset1c{}_4d^3_4L^i_{\underset\vee{\alpha
        n}}a^\alpha_{j\underset3|m}
        +2\overset1c{}_5d^1_5L^i_{\underset\vee{\alpha m}}a^\alpha_{j\underset1|n}
        +2\overset1c{}_5d^2_5L^i_{\underset\vee{\alpha m}}a^\alpha_{j\underset2|n}
        +2\overset1c{}_5d^3_5L^i_{\underset\vee{\alpha
        m}}a^\alpha_{j\underset3|n}
          \\&
          +a^\alpha_j\Big[R^i_{\alpha mn}
          +\overset1c{}_6L^i_{\underset\vee{\alpha m}|n}
          +\overset1c{}_7L^i_{\underset\vee{\alpha n}|m}
          +\big\{\overset1c{}_8-\overset1c{}_4(d^1_4-d^2_4+d^3_4)\big\}
          L^\beta_{\underset\vee{\alpha
          m}}L^i_{\underset\vee{\beta n}}
          \\&+\big\{\overset1c{}_9-
          \overset1c{}_5(d^1_5-d^2_5+d^3_5)\big\}L^\beta_{\underset\vee{\alpha
          n}}L^i_{\underset\vee{\beta m}}
          +\big\{2\overset1c{}_{10}-\overset1c{}_3(d^1_3+d^2_3-d^3_3)\big\}L^i_{\underset\vee{\alpha
          \beta}}L^\beta_{\underset\vee{mn}}
          \Big]\\&
          -a^i_\alpha
          \Big[
          R^\alpha_{jmn}+\overset1c{}_{11}L^\alpha_{\underset\vee{jm}|n}
          +\overset1c{}_{12}L^\alpha_{\underset\vee{jn}|m}
          +\big\{\overset1c{}_{13}-\overset1c{}_1
          (d^1_1-d^2_1-d^3_1)\big\}L^\alpha_{\underset\vee{\beta
          n}}L^\beta_{\underset\vee{jm}}
          \\&
          +\big\{\overset1c{}_{14}-\overset1c{}_2(d^1_2-d^2_2-d^3_2)\big\}L^\alpha_{\underset\vee{\beta
          m}}L^\beta_{\underset\vee{jn}}
          +\big\{2\overset1c{}_{15}-\overset1c{}_3
          (d^1_3-d^2_3-d^3_3)\big\}L^\alpha_{\underset\vee{j\beta}}L^\beta_{\underset\vee{mn}}
          \Big]\\&
          -2a^\alpha_\beta\Big[
          \big\{\overset1c{}_{16}+\overset1c{}_2(d^1_2-d^2_2+d^3_2)
          -\overset1c{}_5(d^1_5-d^2_5-d^3_5)\big\}L^i_{\underset\vee{\alpha
          m}}L^\beta_{\underset\vee{jn}}
          \\&
          +\big\{\overset1c{}_{17}+
          \overset1c{}_1(d^1_1-d^2_1+d^3_1)
          -\overset1c{}_4(d^1_4-d^2_4-d^3_4)\big\}L^i_{\underset\vee{\alpha
          n}}L^\beta_{\underset\vee{jm}}
          \Big],
        \endaligned\label{eq:riccitypeid123}
        \end{equation}
        \end{scriptsize}

        \noindent  of the identities of Ricci Type, for
        real constants $d^1_k,d^2_k,d^3_k$, $d^1_k+d^2_k+d^3_k=1,d^l_k\in\mathbb R,
        k=1,\ldots,5,l=1,2,3$, are equivalent.\qed
  \end{cor}

  With respect to the definition of the curvature tensor
  for the space $\mathbb{GA}_N$ \big($R(X;Y,Z)=\nabla_Z\nabla_YX-\nabla_Y\nabla_ZX+
  \nabla_{[Y,Z]}X$\big), we conclude that the right sides of the
  identities of Ricci Type (\ref{eq:vricciti}) are the curvatures
  for the space $\mathbb{GA}_N$ expressed in the standard base.

  \begin{rem}
    The covariant derivatives \emph{(\ref{eq:covderivativesim}, \ref{eq:covderivativensim})} are expressed as vectors. We proved that three of these five vectors are linearly independent. For this reason, the double covariant derivatives \emph{(\ref{eq:||11}---\ref{eq:||33})} are vectors. Finally, the right sides of the corresponding identities of Ricci Type are treated as vectors. For this reason, the left and right sides of the Ricci Type identities are obtained as the linear combinations of the corresponding sides in the equations \emph{(\ref{eq:ric11-11}---\ref{eq:ric33-33})}.
    Hence, when we say that some Ricci Type identities are linearly (in)dependent
    we think that for the corresponding equalities.

    Note also that the equation \emph{(\ref{eq:riccitypeid123})} holds with respect to the equalities

    \begin{eqnarray*}
      a^i_{j\underset1|k}=a^i_{j|k}+
      L^i_{\underset\vee{\alpha k}}a^\alpha_j-
      L^\alpha_{\underset\vee{jk}}a^i_\alpha,&
      a^i_{j\underset2|k}=a^i_{j|k}-
      L^i_{\underset\vee{\alpha k}}a^\alpha_j+
      L^\alpha_{\underset\vee{jk}}a^i_\alpha,&
      a^i_{j\underset3|k}=a^i_{j|k}+
      L^i_{\underset\vee{\alpha k}}a^\alpha_j+
      L^\alpha_{\underset\vee{jk}}a^i_\alpha,
    \end{eqnarray*}

    \noindent substituted into the equation
    \emph{(\ref{eq:vricciti})}. In the same manner, we may change
    the covariant derivatives $L^i_{\underset\vee{jm}|n}$ to
    linear combinations of the covariant derivatives
    $L^i_{\underset\vee{jm}\underset1|n}$,
    $L^i_{\underset\vee{jm}\underset2|n}$,
    $L^i_{\underset\vee{jm}\underset3|n}$, but we preserved the
    covariant derivative $L^i_{\underset\vee{jm}|n}$ in the
    equations \emph{(\ref{eq:vricciti}, \ref{eq:riccitypeid123})}
    with respect to the results from the previously published papers
    about linearly independent curvature tensors. Because linear
    combinations of tensors are tensors, such as the linear combinations
    of the covariant derivatives $L^i_{\underset\vee{jm}|n}$,
    \ldots, $L^i_{\underset\vee{jm}\underset4|n}$ as well, this transformation would
    not change the tensor characteristics of the curvature tensors for the space
    $\mathbb{GA}_N$. The curvature pseudotensors are lost through
    the computation and they composed by the tensor $a^i_j$ are not
    components of the curvatures for the space $\mathbb{GA}_N$.
  \end{rem}

  In the rows above, we obtained that curvature tensors for the
  non-symmetric affine connection space $\mathbb{GA}_N$ together with its
  torsion tensor are enough to express the family of curvatures for this space.

  We will discuss the curvature characteristics
  of pseudocurvature tensors obtained in \cite{cvetkoviczlatanovic1, mincic2, mincic4,
  mincvel, mincicnovi, z4} in the subsection
  \ref{subsectiontransformationctpct}.

  \subsection{Linearly independent curvatures}

  From the linearly independent Ricci-Type identities, it is
  obtained fourteen curvature tensors. They are given by the
  equations
  (\ref{eq:rho1}---\ref{eq:rho14}) from the Appendix II. As in the previous works
  \cite{cvetkoviczlatanovic1, mincic2, mincic4, mincicnovi, mincvel, z4}
  and many others, these
  tensors are elements of the family

  \begin{equation}
    \overset1\rho{}^i_{jmn}=R^i_{jmn}
    +uL^i_{\underset\vee{jm}|n}
    +u'L^i_{\underset\vee{jn}|m}
    +vL^\alpha_{\underset\vee{jm}}L^i_{\underset\vee{\alpha n}}
    +v'L^\alpha_{\underset\vee{jn}}L^i_{\underset\vee{\alpha m}}
    +wL^\alpha_{\underset\vee{mn}}L^i_{\underset\vee{j\alpha}},
    \label{eq:rhofamily}
  \end{equation}

  \noindent for the corresponding coefficients $u,u',v,v',w$.

  It holds the following theorem and its corollary.
  \begin{thm}
    Six of the curvature tensors from the family
    \emph{(\ref{eq:rhofamily})} are linearly independent. \qed
  \end{thm}

  \begin{cor}
    The curvature tensors $\underset1{\overset1\rho}{}^i_{jmn}$,
    $\overset1{\underset2\rho}{}^i_{jmn}$,
    $\overset1{\underset3\rho}{}^i_{jmn}$,
    $\overset1{\underset4\rho}{}^i_{jmn}$,
    $\overset1{\underset7\rho}{}^i_{jmn}$,
    $\overset1{\underset{10}\rho}{}^i_{jmn}$, given by the equations
    \emph{(\ref{eq:rho1})}, \emph{(\ref{eq:rho2})},
    \emph{(\ref{eq:rho3})}, \emph{(\ref{eq:rho4})},
    \emph{(\ref{eq:rho7})}, \emph{(\ref{eq:rho10})}, in the Appendix II, are linearly
    independent.

    The curvature tensors $R^i_{jmn}$,
    $\overset1{\underset2\rho}{}^i_{jmn}$,
    $\overset1{\underset3\rho}{}^i_{jmn}$,
    $\overset1{\underset4\rho}{}^i_{jmn}$,
    $\overset1{\underset7\rho}{}^i_{jmn}$,
    $\overset1{\underset{10}\rho}{}^i_{jmn}$, given by
    the equation \emph{(\ref{eq:R})} and the equations
    \emph{(\ref{eq:rho2})},
    \emph{(\ref{eq:rho3})},
    \emph{(\ref{eq:rho4})},
    \emph{(\ref{eq:rho7})},
    \emph{(\ref{eq:rho10})}, in the Appendix II, are linearly
    independent.

    The curvature tensors $\underset1{\overset1\rho}{}^i_{jmn}$,
    $\overset1{\underset2\rho}{}^i_{jmn}$,
    $\overset1{\underset3\rho}{}^i_{jmn}$,
    $\overset1{\underset4\rho}{}^i_{jmn}$,
    $\overset1{\underset7\rho}{}^i_{jmn}$,
    $R^i_{jmn}$, given by the equations
    \emph{(\ref{eq:rho1})},
    \emph{(\ref{eq:rho2})},
    \emph{(\ref{eq:rho3})},
    \emph{(\ref{eq:rho4})},
    \emph{(\ref{eq:rho7})}, in the Appendix II,
    and the equation \emph{(\ref{eq:R})}, are linearly
    independent.\qed
  \end{cor}
  \subsection{Curvature pseudotensors}\label{subsectiontransformationctpct}

  It is obtained in \cite{mincic2, mincic4}, and confirmed in many papers
  after \big(see for example \cite{cvetkoviczlatanovic1, mincicnovi, mincvel, z4}\big),
  the existence of the fifteen curvature pseudotensors for the space $\mathbb{GA}_N$.

  The main characteristic of the fourteen identities of Ricci Type where the
  curvature pseudotensors are obtained is that they depend of
  partial derivatives of the tensor $a^i_j$.

  If we substitute the definition (\ref{eq:covderivativesim}) of the
  covariant derivative with respect to symmetric affine connection into
  the equation (\ref{eq:vricciti}), we will obtain the equation
\begin{scriptsize}
  \begin{equation}
    \aligned
        a^i_{j\underset p|m\underset q|n}-
        a^i_{j\underset r|n\underset s|m}&=
        \overset 1c{}_1L^\alpha_{\underset\vee{jm}}a^i_{\alpha,n}
        +\overset 1c{}_2L^\alpha_{\underset\vee{jn}}a^i_{\alpha,m}
        +\overset 1c{}_3L^\alpha_{\underset\vee{mn}}a^i_{j,\alpha}
        +\overset 1c{}_4L^i_{\underset\vee{\alpha n}}a^\alpha_{j,m}
        +\overset 1c{}_5L^i_{\underset\vee{\alpha
        m}}a^\alpha_{j,n}\\
        &+a^\alpha_j\big(
        R^i_{\alpha mn}+\overset1c{}_6L^i_{\underset\vee{\alpha
        m}|n}+\overset 1c{}_7L^i_{\underset\vee{\alpha n}|m}
        +\overset 1c{}_8L^\beta_{\underset\vee{\alpha
        m}}L^i_{\underset\vee{\beta n}}+\overset
        1c{}_9L^\beta_{\underset\vee{\alpha
        n}}L^i_{\underset\vee{\beta
        m}}+2\overset1c{}_{10}L^i_{\underset\vee{\alpha\beta}}L^\beta_{\underset\vee{mn}}\\&
        +\overset1c{}_3L^\beta_{\underset\vee{mn}}L^i_{\underline{\alpha\beta}}
        +\overset1c{}_4L^i_{\underset\vee{\beta
        n}}L^\beta_{\underline{\alpha m}}
        +\overset1c{}_5L^i_{\underset\vee{\beta
        m}}L^\beta_{\underline{\alpha n}}
        \big)\\&
        -a^i_\alpha\big(R^\alpha_{jmn}+\overset1c{}_{11}L^\alpha_{\underset\vee{jm}|n}
        +\overset1c{}_{12}L^\alpha_{\underset\vee{jn}|m}
        +\overset1c{}_{13}L^\alpha_{\underset\vee{\beta
        n}}L^\beta_{\underset\vee{jm}}
        +\overset1c{}_{14}L^\alpha_{\underset\vee{\beta
        m}}L^\beta_{\underset\vee{jn}}
        +\overset1c{}_{15}L^\alpha_{\underset\vee{j\beta
        }}L^\beta_{\underset\vee{mn}}\\&
        +\overset1c{}_1L^\beta_{\underset\vee{jm}}L^\alpha_{\underline{\beta
        n}}+\overset1c{}_2L^\beta_{\underset\vee{jn}}L^\alpha_{\underline{\beta
        m}}+\overset1c{}_3L^\beta_{\underset\vee{mn}}L^\alpha_{\underline{j\beta
        }}
        \big)\\&
        -a^\alpha_\beta\big(\overset1c{}_{16}L^i_{\underset\vee{\alpha
        m}}L^\beta_{\underset\vee{jn}}+
        \overset1c{}_{17}L^i_{\underset\vee{\alpha
        n}}L^\beta_{\underset\vee{jm}}
        \!-\!\overset1c{}_1L^\beta_{\underset\vee{jm}}L^i_{\underline{\alpha
        n}}\!-\!
        \overset1c{}_2L^\beta_{\underset\vee{jn}}L^i_{\underline{\alpha
        m}}\!+\!
        \overset1c{}_4L^i_{\underset\vee{\beta n}}L^\alpha_{\underline{jm
        }}\!+\!
        \overset1c{}_5L^i_{\underset\vee{\beta
        m}}L^\alpha_{\underline{jn
        }}
        \big).
        \endaligned\label{eq:vvricciti}
  \end{equation}
  \end{scriptsize}

  Let us consider the geometrical objects $a^i_{j<mn>}$,
  $a^i_{j\leqslant mn\geqslant}$, $a^i_{j\eqslantless
  mn\eqslantgtr}$, $a^i_{j\leqslant mn\eqslantgtr}$,
  $a^i_{j\eqslantless mn\geqslant}$,
   defined in \cite{mincic4,
  mincic2} and used in later papers. We also need to
  use the equality $L^i_{jk}=L^i_{\underline{jk}}+
  L^i_{\underset\vee{jk}}$ in the following computations.

  Hence, one gets
  \begin{align}
    &a^i_{j<mn>}=L^i_{\underset\vee{\alpha m}}a^\alpha_{j,n}-
    L^\alpha_{\underset\vee{jm}}a^i_{\alpha,n},\label{eq:<>1}\\
    &a^i_{j\leqslant mn\geqslant}=\big(L^i_{m\alpha}L^\beta_{jn}-
    L^i_{\alpha m}L^\beta_{nj}\big)a^\alpha_\beta=
    \big(L^i_{\underline{\alpha m}}L^\beta_{\underset\vee{jn}}
    -L^i_{\underset\vee{\alpha
    m}}L^\beta_{\underline{jn}}\big)a^\alpha_\beta,\label{eq:<>2}\\
    &a^i_{j\eqslantless mn\eqslantgtr}=
    \big(L^i_{m\alpha}L^\beta_{nj}-L^i_{\alpha
    m}L^\beta_{jn}\big)a^\alpha_\beta=-\big(L^i_{\underline{\alpha
    m}}L^\beta_{\underset\vee{jn}}+L^i_{\underset\vee{\alpha
    m}}L^\beta_{\underline{jn}}\big)a^\alpha_\beta,\label{eq:<>3}\\
    &
    a^i_{j\leqslant mn\eqslantgtr}=
    \big(L^i_{\underset\vee{m\alpha}}L^\beta_{jn}\!-\!
    L^i_{\alpha n}L^\beta_{\underset\vee{mj}}\big)a^\alpha_\beta=
    -\big(L^i_{\underset\vee{\alpha
    m}}L^\beta_{\underset\vee{jn}}\!-\!
    L^i_{\underset\vee{\alpha
    n}}L^\beta_{\underset\vee{jm}}\big)a^\alpha_\beta
    \!-\!\big(L^i_{\underset\vee{\alpha
    m}}L^\beta_{\underline{jn}}\!-\!
    L^i_{\underline{\alpha
    n}}L^\beta_{\underset\vee{jm}}\big)a^\alpha_\beta,\label{eq:<>4}\\
    &
    a^i_{j\eqslantless
    mn\geqslant}=\big(L^i_{m\alpha}L^\beta_{\underset\vee{jn}}-L^i_{\underset\vee{\alpha
    n}}L^\beta_{mj}\big)a^\alpha_\beta=
    -\big(L^i_{\underset\vee{\alpha
    m}}L^\beta_{\underset\vee{jn}}\!-\!
    L^i_{\underset\vee{\alpha
    n}}L^\beta_{\underset\vee{jm}}\big)a^\alpha_\beta\!+\!
    \big(L^i_{\underline{\alpha m}}L^\beta_{\underset\vee{jn}}\!-\!
    L^i_{\underset\vee{\alpha
    n}}L^\beta_{\underline{jm}}\big)a^\alpha_\beta,\label{eq:<>5}\\
    &\aligned
    L^\alpha_{mn}a^i_{j\underset1|\alpha}
    &\overset{(\ref{eq:covderivativensim})}=L^\alpha_{\underset\vee{mn}}
    \big(a^i_{j,\alpha}+L^i_{\underline{\beta\alpha}}a^\beta_j-
    L^\beta_{\underline{j\alpha}}a^i_\beta\big)+
    L^\alpha_{\underset\vee{mn}}\big(L^i_{\underset\vee{\beta\alpha}}a^\beta_j-
    L^\beta_{\underset\vee{j\alpha}}a^i_\beta\big)\\&
    +L^\alpha_{\underline{mn}}
    \big(a^i_{j,\alpha}+L^i_{\underline{\beta\alpha}}a^\beta_j-
    L^\beta_{\underline{j\alpha}}a^i_\beta\big)+
    L^\alpha_{\underline{mn}}\big(L^i_{\underset\vee{\beta\alpha}}a^\beta_j-
    L^\beta_{\underset\vee{j\alpha}}a^i_\beta\big),
    \endaligned\label{eq:<>6}\\
    &\aligned
    L^\alpha_{mn}a^i_{j\underset2|\alpha}
    &\overset{(\ref{eq:covderivativensim})}=L^\alpha_{\underset\vee{mn}}
    \big(a^i_{j,\alpha}+L^i_{\underline{\alpha\beta}}a^\beta_j-
    L^\beta_{\underline{\alpha j}}a^i_\beta\big)+
    L^\alpha_{\underset\vee{mn}}\big(L^i_{\underset\vee{\alpha\beta}}a^\beta_j-
    L^\beta_{\underset\vee{\alpha j}}a^i_\beta\big)\\&
    +L^\alpha_{\underline{mn}}
    \big(a^i_{j,\alpha}+L^i_{\underline{\alpha\beta}}a^\beta_j-
    L^\beta_{\underline{\alpha j}}a^i_\beta\big)+
    L^\alpha_{\underline{mn}}\big(L^i_{\underset\vee{\alpha\beta}}a^\beta_j-
    L^\beta_{\underset\vee{\alpha j}}a^i_\beta\big).
    \endaligned\label{eq:<>7}
  \end{align}

  After recognizing the geometrical objects
  (\ref{eq:<>1}---\ref{eq:<>5}) in the equation (\ref{eq:vvricciti})
  and adding the left sides of the equations (\ref{eq:<>6},
  \ref{eq:<>7}) but subtracting the right sides of them from the
  brackets in the second, third, fourth and fifth row of this
  equation, one obtains all the curvature pseudotensors searched in the
  papers \cite{mincic2,
  mincic4,cvetkoviczlatanovic1,mincicnovi,mincvel,z4}.

  The curvatures of the space $\mathbb{GA}_N$ are the right sides of
  the equations (\ref{eq:vricciti}), for different $p,q,r,s$. With
  respect to these equations, we are able to obtain the curvature $a^i_{j|m|n}-a^i_{j|n|m}$
   of the associated space $\mathbb A_N$. The curvature tensor
   $R^i_{jmn}$ of the space $\mathbb A_N$ uniquely defines its
   curvature. Any of summands
   $\overset1c{}_dL^i_{\underset\vee{jk}}a^u_{v|c}$,
   $d=1,\ldots,5$ from the equation (\ref{eq:vricciti}) may be obtained as the corresponding linear
   combination of the differences $a^i_{j\underset p|m\underset
   q|n}-a^i_{j\underset r|n\underset s|m}$ and the curvature $R^\alpha_{jmn}a^i_\alpha-
   R^i_{\alpha mn}a^\alpha_j$. For this reason, any of them is
   the component of curvature for the space $\mathbb{GA}_N$.

   No one of the summands
   $L^i_{\underset\vee{jk}}L^p_{\underline{qr}}a^u_v$ which generate
   curvature pseudotensors cannot be obtained from linear
   combinations of the differences $a^i_{j\underset p|m\underset
   q|n}-a^i_{j\underset r|n\underset s|m}$ separately from the
   covariant derivatives. That means that these summands indirectly
   generate the curvatures of the space $\mathbb{GA}_N$ (together
   with the partial derivatives $a^i_{j,k}$, they are
   the components of the curvatures
    $\overset1c{}_dL^i_{\underset\vee{jk}}a^u_{v|c}$). Hence, the
   curvature tensors together with the torsion tensor uniquely
   determine the curvature of the space $\mathbb{GA}_N$. The
   curvature pseudotensors do not have the same role as the
   curvature tensors for the space $\mathbb{GA}_N$.

   The same statement holds for the family (\ref{eq:riccitypeid123}).

   To point, curvature pseudotensors of the space $\mathbb{GA}_N$
   may anticipate linear combinations of components of curvatures for this space.
   However, they are linear combinations of one complete and one
   incomplete component of this curvature. For this reason, we need
   to solve a system of equations to find the curvature with respect
   to curvature pseudotensors.

    \section{EXAMPLE: Application of curvature tensors of space $\mathbb{GA}_N$}

    Matthias Blau recalled different published results to gather different
    findings in the field of cosmology \big(see \cite{blau}\big).
    From the other side, I. Shapiro \cite{shapiro} and many other authors have studied cosmology
    with respect to torsion.
    We will correlate these results below.

    Let us consider the generalized Riemannian spacetime $\mathbb{GR}_4$
    equipped by the metric

    \begin{equation}
      \big(b_{ij}\big)=\left[\begin{array}{cccc}
        s_1(t)&0&0&0\\
        0&s_2(t)&n(t)&0\\
        0&-n(t)&s_3(t)&0\\
        0&0&0&s_4(t)
      \end{array}\right],
      \label{eq:b1metrictorsion}
    \end{equation}

    \noindent for the differentiable functions
    $s_1(t),\ldots,s_4(t),n(t)$.

    The symmetric and anti-symmetric parts of the metric tensor
    $b_{ij}$ are

    \begin{eqnarray}
      \big(b_{\underline{ij}}\big)=\left[\begin{array}{cccc}
        s_1(t)&0&0&0\\
        0&s_2(t)&0&0\\
        0&0&s_3(t)&0\\
        0&0&0&s_4(t)
      \end{array}\right]
      &\mbox{and}&
      \big(b_{\underset\vee{ij}}\big)=\left[\begin{array}{cccc}
        0&0&0&0\\
        0&0&n(t)&0\\
        0&-n(t)&0&0\\
        0&0&0&0
      \end{array}\right].
      \label{eq:bsimantisim}
    \end{eqnarray}

    The covariant anti-symmetric part of the corresponding
    generalized Christoffel symbol is

    \begin{equation}
      \Gamma_{1.\underset\vee{23}}=
      -\Gamma_{1.\underset\vee{32}}=
      -\Gamma_{2.\underset\vee{13}}=
      \Gamma_{2.\underset\vee{31}}=
      \Gamma_{3.\underset\vee{12}}=
      -\Gamma_{3.\underset\vee{21}}=-\frac12n'(t),
      \label{eq:gammabantisim}
    \end{equation}

    \noindent and $\Gamma_{i.\underset\vee{jk}}=0$ in all other
    cases.

    \subsection{Motivation from cosmology}

    In the Blau's book \cite{blau}, it is analyzed the
    Einstein-Hilbert action

    \begin{equation}
      S=\int{d^4x\sqrt{|b|}\big(R+\mathcal L_M\big)},
      \label{eq:ehactionblau}
    \end{equation}

    \noindent for the scalar
    curvature
    $R=g^{\underline{\alpha\beta}}R^\gamma_{\alpha\beta\gamma}$ of
    the Riemannian space $\mathbb R_4$, the metric determinant
    $b=\det\big(b_{\underline{ij}}\big)$ and
    part $\mathcal L_M$ describing
    any matter fields appearing in the theory.

    With respect to the equation (\ref{eq:ehactionblau}), it is
    obtained the Einstein's equations of motion
    \begin{equation}
    R_{ij}-\frac12Rb_{\underline{ij}}=T_{ij},\label{eq:einsteinmotionblau}
    \end{equation}

    \noindent for the
    energy-momentum tensor $T_{ij}$.
    The Energy-Momentum Tensors $T_{ij}$
    with respect to different parts $\mathcal L_M$ are recalled in
    \cite{blau}.

    From the other side, I. Shapiro \cite{shapiro} have studied the
    cosmology with respect to torsion. The results in Shapiro's
    article are equivalent to the results from \cite{blau}. The
    Shapiro's
    findings are correlated with torsion unlike the results in
    \cite{blau}.

    Our purpose is to recognize some of the expressions
    from the Blau's and Shapiro's works with
    respect to the curvature tensors of the
     above defined generalized spacetime  $\mathbb{GR}_4$
    in this section.

    \subsection{Physical interpretation of curvature tensors and torsion}

    The family $\overset1\rho=b^{\underline{\alpha\beta}}
    \overset1\rho{}^\gamma_{\alpha\beta\gamma}$ of scalar curvatures
    obtained from the curvature tensors
    (\ref{eq:rho1}---\ref{eq:rho14}), Appendix II,
    for the space $\mathbb{GR}_4$ is

    \begin{equation}
      \overset1\rho=R+(v'-w)b^{\underline{\alpha\beta}}
      b^{\underline{\gamma\epsilon}}
      b^{\underline{\delta\zeta}}
      \Gamma_{\alpha.\underset\vee{\gamma\delta}}
      \Gamma_{\beta.\underset\vee{\epsilon\zeta}},
      \label{eq:physscalarcurvaturefamily}
    \end{equation}

    \noindent with respect to the corresponding coefficients $v'$
    and $w$.

    Based on the considerations from the Shapiro's work
    \big(\cite{shapiro}, section 2\big), we conclude that the family of the Einstein-Hilbert actions
    with torsion is

    \begin{equation}
      \widetilde S=\int{d^4x\sqrt{|b|}\big(R+
      (v'-w)b^{\underline{\alpha\beta}}
      b^{\underline{\gamma\epsilon}}
      b^{\underline{\delta\zeta}}
      \Gamma_{\alpha.\underset\vee{\gamma\delta}}
      \Gamma_{\beta.\underset\vee{\epsilon\zeta}}\big)}.
      \label{eq:ehactiontorsion}
    \end{equation}

    After comparing the equations (\ref{eq:ehactionblau},
    \ref{eq:ehactiontorsion}), one obtains

    \begin{equation}
      \mathcal L_M=(v'-w)b^{\underline{\alpha\beta}}
      b^{\underline{\gamma\epsilon}}
      b^{\underline{\delta\zeta}}
      \Gamma_{\alpha.\underset\vee{\gamma\delta}}
      \Gamma_{\beta.\underset\vee{\epsilon\zeta}}.
      \label{eq:LMtorsion}
    \end{equation}

    For the further work, we need to recall the term of \emph{the functional derivative}
    \cite{greiner}.
    \begin{enumerate}[-]
      \item Consider the functional

      \begin{equation}
        J[f]=\int_a^b{L\big[f(t),f'(t)\big]dt},
      \end{equation}

      \noindent for $f'(t)=df(t)/dt$. The variational derivative of
      the functional $J[f]$ by the function $f(t)$ is

      \begin{equation}
        \frac{\delta J[f]}{\delta
        f(t)}=\int{dt\Big(\frac{\partial L\big[f(t),f'(t)\big]}{\partial f(t)}-
        \dfrac d{dt}\frac{\partial
        L\big[f(t),f'(t)\big]}{f'(t)}\Big)}
        =\lim_{\varepsilon\to0}{\frac{J\big[f(t)+\varepsilon\delta(x-t)\big]-J[f]}\varepsilon}
        ,
        \label{eq:funderivative}
      \end{equation}

      \noindent for the Dirac $\delta$-function $\delta(x)$ and $x\neq t$. As in
      \cite{greiner}, the limit $\varepsilon\to0$ has to be taken
      first, before other limiting operations.

      The variation of the functional $J[f]$ is
      \begin{equation}
        \delta J[f]=\int{\frac{\delta J[f]}{\delta f(t)}\delta
        f(t)dt}.\label{eq:Jvariation}
      \end{equation}
      \item Note the following equalities
      \begin{align}
        &\frac{\delta\Big\{F[f]\cdot G[f]\Big\}}{\delta f(t)}=
        \frac{\delta F[f]}{\delta f(t)}\cdot G[f]+
        \frac{\delta G[f]}{\delta f(t)}\cdot F[f],\label{funderivativeproperties1}
        \\\displaybreak[0]&
        \frac{\delta}{\delta
        f(t)}F\big[G[f]\big]=\int{dt\frac{\delta F[G]}{\delta
        G(t)}\cdot \frac{\delta{G[f]}}{\delta f(t)}},
        \label{eq:funderivativeproperties2}\\\displaybreak[0]
        &\frac{\delta F[f]}{\delta
        C}=0,\label{eq:funderivativeproperties3}
      \end{align}

      for a constant function $f(t)=C$.
    \end{enumerate}

    Based on the equations (\ref{eq:bsimantisim}, \ref{eq:gammabantisim}, \ref{eq:LMtorsion}), we get

    \begin{equation}
        \mathcal L_M=\frac{3(v'-w)}2\cdot
        \frac{\big(n'(t)\big)^2}{s_1(t)s_2(t)s_3(t)}.
        \label{eq:LMb}
    \end{equation}

    It is evident that the components

    \begin{eqnarray}
      n_{1,2}(t)=\pm\frac{2}{3(v'-w)}\int{\sqrt{\mathcal L_M\cdot s_1(t) s_2(t)
      s_3(t)}dt},
      \label{eq:2nLM}
    \end{eqnarray}

    \noindent of the anti-symmetric part
    $b_{\underset\vee{ij}}$ for the metric tensor (\ref{eq:b1metrictorsion})
     correspond to the same operator $\mathcal L_M$. Geometrically,
    the operator $\mathcal L_M$ generates two generalized Riemannian
    spaces $\mathbb{GR}_4^+$ and $\mathbb{GR}_4^-$ in the sense of Eisenhart's definitions
    \cite{eis01, eis02}. Moreover, the operator $\mathcal L_M$
    generates two opposite torsion tensors $\underset1T{}^i_{jk}=2\Gamma^i_{\underset\vee{jk}}$ and
    $\underset2T{}^i_{jk}=-2\Gamma^i_{\underset\vee{jk}}=-\underset1T{}^i_{jk}$.

    If $n_1(t)=n_2(t)$ in the equation (\ref{eq:2nLM}), we get $\mathcal L_M=0$.
    Equivalently, the equality $n_1(t)=n_2(t)$
    corresponds to the no matter part of the space. The physical
    considerations about the cases of $\mathcal L_M=0$
    are
     geometrically covered by the corresponding Riemannian spaces $\mathbb
     R_4$.

    After vanishing the variation (\ref{eq:Jvariation}) of the
    Einstein-Hilbert action (\ref{eq:ehactionblau}) and using the
    equation (\ref{eq:LMtorsion}), we
    get the family of Energy-Momentum Tensors $T_{ij}=
    -2\delta\mathcal L_M/\delta
    g^{\underline{ij}}+b_{\underline{ij}}\mathcal L_M$, i.e.

    \begin{equation}
      T_{ij}=
      -3(v'-w)\frac{\delta\Big(\frac{\big(n'(t)\big)^2}{s_1(t)s_2(t)s_3(t)}\Big)}
      {\delta
      b^{\underline{ij}}}+\frac{3(v'-w)}2b_{\underline{ij}}\cdot
      \frac{\big(n'(t)\big)^2}{s_1(t)s_2(t)s_3(t)},
    \end{equation}

    \noindent with respect to the family of scalar curvatures
    (\ref{eq:physscalarcurvaturefamily}) of the space $\mathbb{GR}_4$.

    For researches in physics, it is the most
     common to examine the case of $v'-w=1$ for the coefficients $v'$ and $w$
     in the equation (\ref{eq:physscalarcurvaturefamily}).

    \section{Conclusion}

    We achieved the aims of this paper above.

    In the section 2, we proved that three kinds of
    covariant derivatives are enough to be defined for the complete
    analysis of a non-symmetric affine connection. Moreover, we
    obtained that there are seventeen linearly independent
    Ricci-Type identities and six linearly independent curvature
    tensors of the non-symmetric affine connection space $\mathbb{GA}_N$. In that section, it is
    explained why the curvature pseudotensors are not components of curvatures
    for
    this space.

    In the section 3, we physically interpreted the curvature
    tensors obtained in the section 2. Namely, we founded that
    the anti-symmetric parts of affine connection coefficients
    correspond to matter and obtained the
    corresponding family of Energy-Momentum Tensors. The results
    obtained in that section motivate the author to find the
    general formulae for pressure, energy-density and state
    parameter of a cosmological fluid.

    In the future, we are aimed to generalize the results from the
    previous paper about invariants of geometric mappings. Moreover,
    we will try to apply differential geometry in physics, specially
    in cosmology, more detail than in this paper.
    \section{Appendix I: Linearly independent identities of Ricci Type}
    \begin{align}
      &\aligned
      a^i_{j\underset1|m\underset1|n}-
      a^i_{j\underset1|n\underset1|m}&=-2L^\alpha_{\underset\vee{mn}}a^i_{j|\alpha}
      \\&+a^\alpha_j\big(R^i_{\alpha mn}+
      L^i_{\underset\vee{\alpha m}|n}
      -L^i_{\underset\vee{\alpha n}|m}
      +L^\beta_{\underset\vee{\alpha m}}L^i_{\underset\vee{\beta n}}
      -L^\beta_{\underset\vee{\alpha n}}L^i_{\underset\vee{\beta m}}
      -2L^i_{\underset\vee{\alpha\beta}}L^\beta_{\underset\vee{mn}}\big)\\&
      -a^i_\alpha\big(R^\alpha_{jmn}
      +L^\alpha_{\underset\vee{jm}|n}
      -L^\alpha_{\underset\vee{jn}|m}
      +L^\beta_{\underset\vee{jm}}L^\alpha_{\underset\vee{\beta n}}
      -L^\beta_{\underset\vee{jn}}L^\alpha_{\underset\vee{\beta m}}
      -2L^\alpha_{\underset\vee{j\beta}}L^\beta_{\underset\vee{mn}}\big),
      \endaligned\label{eq:ric11-11}\\\displaybreak[0]
      &\aligned
      a^i_{j\underset1|m\underset2|n}-
      a^i_{j\underset1|n\underset1|m}&=
      2L^\alpha_{\underset\vee{jn}}a^i_{\alpha|m}
      -2L^i_{\underset\vee{\alpha n}}a^\alpha_{j|m}
      \\&
      +a^\alpha_j\big(R^i_{\alpha mn}
      +L^i_{\underset\vee{\alpha m}|n}
      -L^i_{\underset\vee{\alpha n}|m}
      -L^\beta_{\underset\vee{\alpha m}}L^i_{\underset\vee{\beta n}}
      -L^\beta_{\underset\vee{\alpha n}}L^i_{\underset\vee{\beta m}}
      \big)\\&
      -a^i_\alpha\big(R^\alpha_{jmn}
      +L^\alpha_{\underset\vee{jm}|n}
      -L^\alpha_{\underset\vee{jn}|m}
      +L^\beta_{\underset\vee{jm}}L^\alpha_{\underset\vee{\beta n}}
      +L^\beta_{\underset\vee{jn}}L^\alpha_{\underset\vee{\beta m}}
      \big)\\&
      +2a^\alpha_\beta\big(L^i_{\underset\vee{\alpha m}}L^\beta_{\underset\vee{jn}}
      +L^i_{\underset\vee{\alpha n}}L^\beta_{\underset\vee{jm}}\big),
      \endaligned\label{eq:ric12-11}\\\displaybreak[0]
      &\aligned
      a^i_{j\underset1|m\underset3|n}-
      a^i_{j\underset1|n\underset1|m}&=2L^\alpha_{\underset\vee{jn}}a^i_{\alpha|m}
      +a^\alpha_j\big(R^i_{\alpha mn}
      +L^i_{\underset\vee{\alpha m}|n}
      -L^i_{\underset\vee{\alpha n}|m}
      +L^\beta_{\underset\vee{\alpha m}}L^i_{\underset\vee{\beta n}}
      -L^\beta_{\underset\vee{\alpha n}}L^i_{\underset\vee{\beta m}}
      \big)\\&
      -a^i_\alpha\big(R^\alpha_{jmn}
      +L^\alpha_{\underset\vee{jm}|n}
      -L^\alpha_{\underset\vee{jn}|m}
      +L^\beta_{\underset\vee{jm}}L^\alpha_{\underset\vee{\beta n}}
      +L^\beta_{\underset\vee{jn}}L^\alpha_{\underset\vee{\beta m}}
      \big)+2a^\alpha_\beta L^i_{\underset\vee{\alpha m}}L^\beta_{\underset\vee{jn}},
      \endaligned\label{eq:ric13-11}\\\displaybreak[0]&\aligned
      a^i_{j\underset2|m\underset1|n}-
      a^i_{j\underset1|n\underset1|m}&=
      2L^\alpha_{\underset\vee{jm}}a^i_{\alpha|n}
      -2L^\alpha_{\underset\vee{mn}}a^i_{j|\alpha}
      -2L^i_{\underset\vee{\alpha m}}a^\alpha_{j|n}\\&
      +a^\alpha_j\big(R^i_{\alpha mn}
      -L^i_{\underset\vee{\alpha m}|n}
      -L^i_{\underset\vee{\alpha n}|m}
      -L^\beta_{\underset\vee{\alpha m}}L^i_{\underset\vee{\beta n}}
      -L^\beta_{\underset\vee{\alpha n}}L^i_{\underset\vee{\beta m}}
      \big)\\&
      -a^i_\alpha\big(R^\alpha_{jmn}
      -L^\alpha_{\underset\vee{jm}|n}
      -L^\alpha_{\underset\vee{jn}|m}
      +L^\beta_{\underset\vee{jm}}L^\alpha_{\underset\vee{\beta n}}
      +L^\beta_{\underset\vee{jn}}L^\alpha_{\underset\vee{\beta m}}
      \big)\\&
      +2a^\alpha_\beta\big(L^i_{\underset\vee{\alpha m}}L^\beta_{\underset\vee{jn}}
      +L^i_{\underset\vee{\alpha n}}L^\beta_{\underset\vee{jm}}\big),
      \endaligned\label{eq:ric21-11}\\\displaybreak[0]
      &\aligned
      a^i_{j\underset2|m\underset2|n}-
      a^i_{j\underset1|n\underset1|m}&=
      2L^\alpha_{\underset\vee{jm}}a^i_{\alpha|n}
      +2L^\alpha_{\underset\vee{jn}}a^i_{\alpha|m}
      -2L^i_{\underset\vee{\alpha n}}a^\alpha_{j|m}
      -2L^i_{\underset\vee{\alpha m}}a^\alpha_{j|n}
      \\&
      +a^\alpha_j\big(R^i_{\alpha mn}
      -L^i_{\underset\vee{\alpha m}|n}
      -L^i_{\underset\vee{\alpha n}|m}
      +L^\beta_{\underset\vee{\alpha m}}L^i_{\underset\vee{\beta n}}
      -L^\beta_{\underset\vee{\alpha n}}L^i_{\underset\vee{\beta m}}
      -2L^i_{\underset\vee{\alpha\beta}}L^\beta_{\underset\vee{mn}}\big)\\&
      -a^i_\alpha\big(R^\alpha_{jmn}
      -L^\alpha_{\underset\vee{jm}|n}
      -L^\alpha_{\underset\vee{jn}|m}
      +L^\beta_{\underset\vee{jm}}L^\alpha_{\underset\vee{\beta n}}
      -L^\beta_{\underset\vee{jn}}L^\alpha_{\underset\vee{\beta m}}
      -2L^\alpha_{\underset\vee{j\beta}}L^\beta_{\underset\vee{mn}}\big),
      \endaligned\label{eq:ric22-11}\\\displaybreak[0]
      &\aligned
      a^i_{j\underset2|m\underset3|n}-
      a^i_{j\underset1|n\underset1|m}&=
      2L^\alpha_{\underset\vee{jm}}a^i_{\alpha|n}
      +2L^\alpha_{\underset\vee{jn}}a^i_{\alpha|m}
      -2L^i_{\underset\vee{\alpha m}}a^\alpha_{j|n}
      \\&
      +a^\alpha_j\big(R^i_{\alpha mn}
      -L^i_{\underset\vee{\alpha m}|n}
      -L^i_{\underset\vee{\alpha n}|m}
      -L^\beta_{\underset\vee{\alpha m}}L^i_{\underset\vee{\beta n}}
      -L^\beta_{\underset\vee{\alpha n}}L^i_{\underset\vee{\beta m}}
      -2L^i_{\underset\vee{\alpha\beta}}L^\beta_{\underset\vee{mn}}\big)\\&
      -a^i_\alpha\big(R^\alpha_{jmn}
      -L^\alpha_{\underset\vee{jm}|n}
      -L^\alpha_{\underset\vee{jn}|m}
      +L^\beta_{\underset\vee{jm}}L^\alpha_{\underset\vee{\beta n}}
      -L^\beta_{\underset\vee{jn}}L^\alpha_{\underset\vee{\beta m}}
      -2L^\alpha_{\underset\vee{j\beta}}L^\beta_{\underset\vee{mn}}\big)
      \\&+2a^\alpha_\beta L^\alpha_{\underset\vee{\alpha n}}L^\beta_{\underset\vee{jm}},
      \endaligned\label{eq:ric23-11}\\\displaybreak[0]
      &\aligned
      a^i_{j\underset3|m\underset1|n}-
      a^i_{j\underset1|n\underset1|m}&=
      2L^\alpha_{\underset\vee{jm}}a^i_{\alpha|n}
      -2L^\alpha_{\underset\vee{mn}}a^i_{j|\alpha}
      \\&
      +a^\alpha_j\big(R^i_{\alpha mn}
      +L^i_{\underset\vee{\alpha m}|n}
      -L^i_{\underset\vee{\alpha n}|m}
      +L^\beta_{\underset\vee{\alpha m}}L^i_{\underset\vee{\beta n}}
      -L^\beta_{\underset\vee{\alpha n}}L^i_{\underset\vee{\beta m}}
      -2L^i_{\underset\vee{\alpha\beta}}L^\beta_{\underset\vee{mn}}\big)\\&
      -a^i_\alpha\big(R^\alpha_{jmn}
      -L^\alpha_{\underset\vee{jm}|n}
      -L^\alpha_{\underset\vee{jn}|m}
      +L^\beta_{\underset\vee{jm}}L^\alpha_{\underset\vee{\beta n}}
      +L^\beta_{\underset\vee{jn}}L^\alpha_{\underset\vee{\beta m}}
      \big)
      +2a^\alpha_\beta L^i_{\underset\vee{\alpha n}}L^\beta_{\underset\vee{jm}},
      \endaligned\label{eq:ric31-11}\\\displaybreak[0]
      &\aligned
      a^i_{j\underset3|m\underset2|n}-
      a^i_{j\underset1|n\underset1|m}&=
      2L^\alpha_{\underset\vee{jm}}a^i_{\alpha|n}
      +2L^\alpha_{\underset\vee{jn}}a^i_{\alpha|m}
      -2L^i_{\underset\vee{\alpha n}}a^\alpha_{j|m}\\&
      +a^\alpha_j\big(R^i_{\alpha mn}+
      L^i_{\underset\vee{\alpha m}|n}
      -L^i_{\underset\vee{\alpha n}|m}
      -L^\beta_{\underset\vee{\alpha m}}L^i_{\underset\vee{\beta n}}
      -L^\beta_{\underset\vee{\alpha n}}L^i_{\underset\vee{\beta m}}
      \big)\\&
      -a^i_\alpha\big(R^\alpha_{jmn}
      -L^\alpha_{\underset\vee{jm}|n}
      -L^\alpha_{\underset\vee{jn}|m}
      +L^\beta_{\underset\vee{jm}}L^\alpha_{\underset\vee{\beta n}}
      -L^\beta_{\underset\vee{jn}}L^\alpha_{\underset\vee{\beta m}}
      -2L^\alpha_{\underset\vee{j\beta}}L^\beta_{\underset\vee{mn}}\big)
      \\&+2a^\alpha_\beta L^i_{\underset\vee{\alpha m}}L^\beta_{\underset\vee{jn}},
      \endaligned\label{eq:ric32-11}\\\displaybreak[0]
      &\aligned
      a^i_{j\underset3|m\underset3|n}-
      a^i_{j\underset1|n\underset1|m}&=
      2L^\alpha_{\underset\vee{jm}}a^i_{\alpha|n}
      +2L^\alpha_{\underset\vee{jn}}a^i_{\alpha|m}
      \\&
      +a^\alpha_j\big(R^i_{\alpha mn}
      +L^i_{\underset\vee{\alpha m}|n}
      -L^i_{\underset\vee{\alpha n}|m}
      +L^\beta_{\underset\vee{\alpha m}}L^i_{\underset\vee{\beta n}}
      -L^\beta_{\underset\vee{\alpha n}}L^i_{\underset\vee{\beta m}}
      \big)\\&
      -a^i_\alpha\big(R^\alpha_{jmn}
      -L^\alpha_{\underset\vee{jm}|n}
      -L^\alpha_{\underset\vee{jn}|m}
      +L^\beta_{\underset\vee{jm}}L^\alpha_{\underset\vee{\beta n}}
      -L^\beta_{\underset\vee{jn}}L^\alpha_{\underset\vee{\beta m}}
      -2L^\alpha_{\underset\vee{j\beta}}L^\beta_{\underset\vee{mn}}\big)
      \\&
      +2a^\alpha_\beta\big(
      L^i_{\underset\vee{\alpha m}}L^\beta_{\underset\vee{jn}}
      +L^i_{\underset\vee{\alpha n}}L^\beta_{\underset\vee{jm}}
      \big),
      \endaligned\label{eq:ric33-11}\\\displaybreak[0]
      &\aligned
      a^i_{j\underset1|m\underset2|n}-
      a^i_{j\underset1|n\underset2|m}&=
      -2L^\alpha_{\underset\vee{jm}}a^i_{\alpha|n}
      +2L^\alpha_{\underset\vee{jn}}a^i_{\alpha|m}
      +2L^\alpha_{\underset\vee{mn}}a^i_{j|\alpha}
      -2L^i_{\underset\vee{\alpha n}}a^\alpha_{j|m}
      +2L^i_{\underset\vee{\alpha m}}a^\alpha_{j|n}
      \\&
      +a^\alpha_j\big(R^i_{\alpha mn}
      +L^i_{\underset\vee{\alpha m}|n}
      -L^i_{\underset\vee{\alpha n}|m}
      -L^\beta_{\underset\vee{\alpha m}}L^i_{\underset\vee{\beta n}}
      +L^\beta_{\underset\vee{\alpha n}}L^i_{\underset\vee{\beta m}}
      +2L^i_{\underset\vee{\alpha\beta}}L^\beta_{\underset\vee{mn}}\big)\\&
      -a^i_\alpha\big(R^\alpha_{jmn}
      +L^\alpha_{\underset\vee{jm}|n}
      -L^\alpha_{\underset\vee{jn}|m}
      -L^\beta_{\underset\vee{jm}}L^\alpha_{\underset\vee{\beta n}}
      +L^\beta_{\underset\vee{jn}}L^\alpha_{\underset\vee{\beta m}}
      +2L^\alpha_{\underset\vee{j\beta}}L^\beta_{\underset\vee{mn}}\big),
      \endaligned\label{eq:ric12-12}\\\displaybreak[0]
      &\aligned
      a^i_{j\underset1|m\underset3|n}-
      a^i_{j\underset1|n\underset2|m}&=
      -2L^\alpha_{\underset\vee{jm}}a^i_{\alpha|n}
      +2L^\alpha_{\underset\vee{mn}}a^i_{j|\alpha}
      +2L^i_{\underset\vee{\alpha n}}a^\alpha_{j|m}
      +2L^i_{\underset\vee{\alpha m}}a^\alpha_{j|n}
      \\&
      +a^\alpha_j\big(R^i_{\alpha mn}
      +L^i_{\underset\vee{\alpha m}|n}
      -L^i_{\underset\vee{\alpha n}|m}
      +L^\beta_{\underset\vee{\alpha m}}L^i_{\underset\vee{\beta n}}
      +L^\beta_{\underset\vee{\alpha n}}L^i_{\underset\vee{\beta m}}
      +2L^i_{\underset\vee{\alpha\beta}}L^\beta_{\underset\vee{mn}}\big)\\&
      -a^i_\alpha\big(R^\alpha_{jmn}
      +L^\alpha_{\underset\vee{jm}|n}
      -L^\alpha_{\underset\vee{jn}|m}
      -L^\beta_{\underset\vee{jm}}L^\alpha_{\underset\vee{\beta n}}
      +L^\beta_{\underset\vee{jn}}L^\alpha_{\underset\vee{\beta m}}
      +2L^\alpha_{\underset\vee{j\beta}}L^\beta_{\underset\vee{mn}}\big)
      \\&-2a^\alpha_\beta L^i_{\underset\vee{\alpha n}}L^\beta_{\underset\vee{jm}},
      \endaligned\label{eq:ric13-12}\\\displaybreak[0]
      &\aligned
      a^i_{j\underset1|m\underset3|n}-
      a^i_{j\underset1|n\underset3|m}&=
      -2L^\alpha_{\underset\vee{jm}}a^i_{\alpha|n}
      +2L^\alpha_{\underset\vee{jn}}a^i_{\alpha|m}
      +2L^\alpha_{\underset\vee{mn}}a^i_{j|\alpha}
      \\&
      +a^\alpha_j\big(R^i_{\alpha mn}
      +L^i_{\underset\vee{\alpha m}|n}
      -L^i_{\underset\vee{\alpha n}|m}
      +L^\beta_{\underset\vee{\alpha m}}L^i_{\underset\vee{\beta n}}
      -L^\beta_{\underset\vee{\alpha n}}L^i_{\underset\vee{\beta m}}
      +2L^i_{\underset\vee{\alpha\beta}}L^\beta_{\underset\vee{mn}}\big)\\&
      -a^i_\alpha\big(R^\alpha_{jmn}
      +L^\alpha_{\underset\vee{jm}|n}
      -L^\alpha_{\underset\vee{jn}|m}
      -L^\beta_{\underset\vee{jm}}L^\alpha_{\underset\vee{\beta n}}
      +L^\beta_{\underset\vee{jn}}L^\alpha_{\underset\vee{\beta m}}
      +2L^\alpha_{\underset\vee{j\beta}}L^\beta_{\underset\vee{mn}}\big)
      \\&+2a^\alpha_\beta
      \big(L^i_{\underset\vee{\alpha m}}L^\beta_{\underset\vee{jn}}
      -L^i_{\underset\vee{\alpha n}}L^\beta_{\underset\vee{jm}}
      \big),
      \endaligned\label{eq:ric13-13}\\\displaybreak[0]
      &\aligned
      a^i_{j\underset2|m\underset1|n}-
      a^i_{j\underset2|n\underset1|m}&=
      2L^\alpha_{\underset\vee{jm}}a^i_{\alpha|n}
      -2L^\alpha_{\underset\vee{jn}}a^i_{\alpha|m}
      -2L^\alpha_{\underset\vee{mn}}a^i_{j|\alpha}
      +2L^i_{\underset\vee{\alpha n}}a^\alpha_{j|m}
      -2L^i_{\underset\vee{\alpha m}}a^\alpha_{j|n}
      \\&
      +a^\alpha_j\big(R^i_{\alpha mn}
      -L^i_{\underset\vee{\alpha m}|n}
      +L^i_{\underset\vee{\alpha n}|m}
      -L^\beta_{\underset\vee{\alpha m}}L^i_{\underset\vee{\beta n}}
      +L^\beta_{\underset\vee{\alpha n}}L^i_{\underset\vee{\beta m}}
      +2L^i_{\underset\vee{\alpha\beta}}L^\beta_{\underset\vee{mn}}\big)\\&
      -a^i_\alpha\big(R^\alpha_{jmn}
      -L^\alpha_{\underset\vee{jm}|n}
      +L^\alpha_{\underset\vee{jn}|m}
      -L^\beta_{\underset\vee{jm}}L^\alpha_{\underset\vee{\beta n}}
      +L^\beta_{\underset\vee{jn}}L^\alpha_{\underset\vee{\beta m}}
      +2L^\alpha_{\underset\vee{j\beta}}L^\beta_{\underset\vee{mn}}\big),
      \endaligned\label{eq:ric21-21}\\\displaybreak[0]
      &\aligned
      a^i_{j\underset2|m\underset2|n}-
      a^i_{j\underset2|n\underset2|m}&=
      2L^\alpha_{\underset\vee{mn}}a^i_{j|\alpha}
      \\&
      +a^\alpha_j\big(R^i_{\alpha mn}
      -L^i_{\underset\vee{\alpha m}|n}
      +L^i_{\underset\vee{\alpha n}|m}
      +L^\beta_{\underset\vee{\alpha m}}L^i_{\underset\vee{\beta n}}
      -L^\beta_{\underset\vee{\alpha n}}L^i_{\underset\vee{\beta m}}
      -2L^i_{\underset\vee{\alpha\beta}}L^\beta_{\underset\vee{mn}}\big)\\&
      -a^i_\alpha\big(R^\alpha_{jmn}
      -L^\alpha_{\underset\vee{jm}|n}
      +L^\alpha_{\underset\vee{jn}|m}
      +L^\beta_{\underset\vee{jm}}L^\alpha_{\underset\vee{\beta n}}
      -L^\beta_{\underset\vee{jn}}L^\alpha_{\underset\vee{\beta m}}
      -2L^\alpha_{\underset\vee{j\beta}}L^\beta_{\underset\vee{mn}}\big),
      \endaligned\label{eq:ric22-22}\\\displaybreak[0]
      &\aligned
      a^i_{j\underset2|m\underset3|n}-
      a^i_{j\underset2|n\underset3|m}&=
      2L^\alpha_{\underset\vee{mn}}a^i_{j|\alpha}
      +2L^i_{\underset\vee{\alpha n}}a^\alpha_{j|m}
      -2L^i_{\underset\vee{\alpha m}}a^\alpha_{j|n}
      \\&
      +a^\alpha_j\big(R^i_{\alpha mn}
      -L^i_{\underset\vee{\alpha m}|n}
      +L^i_{\underset\vee{\alpha n}|m}
      -L^\beta_{\underset\vee{\alpha m}}L^i_{\underset\vee{\beta n}}
      +L^\beta_{\underset\vee{\alpha n}}L^i_{\underset\vee{\beta m}}
      -2L^i_{\underset\vee{\alpha\beta}}L^\beta_{\underset\vee{mn}}\big)\\&
      -a^i_\alpha\big(R^\alpha_{jmn}
      -L^\alpha_{\underset\vee{jm}|n}
      +L^\alpha_{\underset\vee{jn}|m}
      +L^\beta_{\underset\vee{jm}}L^\alpha_{\underset\vee{\beta n}}
      -L^\beta_{\underset\vee{jn}}L^\alpha_{\underset\vee{\beta m}}
      -2L^\alpha_{\underset\vee{j\beta}}L^\beta_{\underset\vee{mn}}\big)
      \\&-2a^\alpha_\beta\big(
      L^i_{\underset\vee{\alpha m}}L^\beta_{\underset\vee{jn}}
      -L^i_{\underset\vee{\alpha n}}L^\beta_{\underset\vee{jm}}
      \big),
      \endaligned\label{eq:ric23-23}\\\displaybreak[0]
      &\aligned
      a^i_{j\underset3|m\underset1|n}-
      a^i_{j\underset3|n\underset1|m}&=
      2L^\alpha_{\underset\vee{jm}}a^i_{\alpha|n}
      -2L^\alpha_{\underset\vee{jn}}a^i_{\alpha|m}
      -2L^\alpha_{\underset\vee{mn}}a^i_{j|\alpha}
      \\&
      +a^\alpha_j\big(R^i_{\alpha mn}
      +L^i_{\underset\vee{\alpha m}|n}
      -L^i_{\underset\vee{\alpha n}|m}
      +L^\beta_{\underset\vee{\alpha m}}L^i_{\underset\vee{\beta n}}
      -L^\beta_{\underset\vee{\alpha n}}L^i_{\underset\vee{\beta m}}
      -2L^i_{\underset\vee{\alpha\beta}}L^\beta_{\underset\vee{mn}}\big)\\&
      -a^i_\alpha\big(R^\alpha_{jmn}
      -L^\alpha_{\underset\vee{jm}|n}
      +L^\alpha_{\underset\vee{jn}|m}
      -L^\beta_{\underset\vee{jm}}L^\alpha_{\underset\vee{\beta n}}
      +L^\beta_{\underset\vee{jn}}L^\alpha_{\underset\vee{\beta m}}
      +2L^\alpha_{\underset\vee{j\beta}}L^\beta_{\underset\vee{mn}}\big)
      \\&-2a^\alpha_\beta\big(
      L^i_{\underset\vee{\alpha m}}L^\beta_{\underset\vee{jn}}
      -L^i_{\underset\vee{\alpha n}}L^\beta_{\underset\vee{jm}}
      \big),
      \endaligned\label{eq:ric31-31}\\\displaybreak[0]
      &\aligned
      a^i_{j\underset3|m\underset3|n}-
      a^i_{j\underset3|n\underset3|m}&=
      2L^\alpha_{\underset\vee{mn}}a^i_{j|\alpha}
      \\&
      +a^\alpha_j\big(R^i_{\alpha mn}
      +L^i_{\underset\vee{\alpha m}|n}
      -L^i_{\underset\vee{\alpha n}|m}
      +L^\beta_{\underset\vee{\alpha m}}L^i_{\underset\vee{\beta n}}
      -L^\beta_{\underset\vee{\alpha n}}L^i_{\underset\vee{\beta m}}
      +2L^i_{\underset\vee{\alpha\beta}}L^\beta_{\underset\vee{mn}}\big)\\&
      -a^i_\alpha\big(R^\alpha_{jmn}
      -L^\alpha_{\underset\vee{jm}|n}
      +L^\alpha_{\underset\vee{jn}|m}
      +L^\beta_{\underset\vee{jm}}L^\alpha_{\underset\vee{\beta n}}
      -L^\beta_{\underset\vee{jn}}L^\alpha_{\underset\vee{\beta m}}
      -2L^\alpha_{\underset\vee{j\beta}}L^\beta_{\underset\vee{mn}}\big).
      \endaligned\label{eq:ric33-33}
    \end{align}

    \section{Appendix II: Curvature tensors obtained from linearly independent
    identities of Ricci Type}
    \setcounter{section}{6}
    \setcounter{equation}{0}
    \begin{align}
      &\overset1{\underset1\rho}{}^i_{jmn}=R^i_{jmn}
      +L^i_{\underset\vee{jm}|n}
      -L^i_{\underset\vee{jn}|m}
      +L^\alpha_{\underset\vee{jm}}L^i_{\underset\vee{\alpha n}}
      -L^\alpha_{\underset\vee{jn}}L^i_{\underset\vee{\alpha m}}
      -2L^i_{\underset\vee{j\alpha}}L^\alpha_{\underset\vee{mn}},\label{eq:rho1}
      \\\displaybreak[0]
      &\overset1{\underset2\rho}{}^i_{jmn}=R^i_{jmn}
      +L^i_{\underset\vee{jm}|n}
      -L^i_{\underset\vee{jn}|m}
      -L^\alpha_{\underset\vee{jm}}L^i_{\underset\vee{\alpha n}}
      -L^\alpha_{\underset\vee{jn}}L^i_{\underset\vee{\alpha m}}
      ,\label{eq:rho2}
      \\\displaybreak[0]
      &\overset1{\underset3\rho}{}^i_{jmn}=R^i_{jmn}
      +L^i_{\underset\vee{jm}|n}
      -L^i_{\underset\vee{jn}|m}
      +L^\alpha_{\underset\vee{jm}}L^i_{\underset\vee{\alpha n}}
      -L^\alpha_{\underset\vee{jn}}L^i_{\underset\vee{\alpha m}}
      ,\label{eq:rho3}
      \\\displaybreak[0]
      &\overset1{\underset4\rho}{}^i_{jmn}=R^i_{jmn}
      -L^i_{\underset\vee{jm}|n}
      -L^i_{\underset\vee{jn}|m}
      -L^\alpha_{\underset\vee{jm}}L^i_{\underset\vee{\alpha n}}
      -L^\alpha_{\underset\vee{jn}}L^i_{\underset\vee{\alpha m}}
      ,\label{eq:rho4}
      \\\displaybreak[0]
      &\overset1{\underset5\rho}{}^i_{jmn}=R^i_{jmn}
      -L^i_{\underset\vee{jm}|n}
      -L^i_{\underset\vee{jn}|m}
      +L^\alpha_{\underset\vee{jm}}L^i_{\underset\vee{\alpha n}}
      -L^\alpha_{\underset\vee{jn}}L^i_{\underset\vee{\alpha m}}
      -2L^i_{\underset\vee{j\alpha}}L^\alpha_{\underset\vee{mn}},\label{eq:rho5}
      \\\displaybreak[0]
      &\overset1{\underset6\rho}{}^i_{jmn}=R^i_{jmn}
      -L^i_{\underset\vee{jm}|n}
      -L^i_{\underset\vee{jn}|m}
      -L^\alpha_{\underset\vee{jm}}L^i_{\underset\vee{\alpha n}}
      -L^\alpha_{\underset\vee{jn}}L^i_{\underset\vee{\alpha m}}
      -2L^i_{\underset\vee{j\alpha}}L^\alpha_{\underset\vee{mn}},\label{eq:rho6}
      \\\displaybreak[0]
      &\overset1{\underset7\rho}{}^i_{jmn}=R^i_{jmn}
      +L^i_{\underset\vee{jm}|n}
      -L^i_{\underset\vee{jn}|m}
      -L^\alpha_{\underset\vee{jm}}L^i_{\underset\vee{\alpha n}}
      +L^\alpha_{\underset\vee{jn}}L^i_{\underset\vee{\alpha m}}
      +2L^i_{\underset\vee{j\alpha}}L^\alpha_{\underset\vee{mn}},\label{eq:rho7}
      \\\displaybreak[0]
      &\overset1{\underset8\rho}{}^i_{jmn}=R^i_{jmn}
      +L^i_{\underset\vee{jm}|n}
      -L^i_{\underset\vee{jn}|m}
      +L^\alpha_{\underset\vee{jm}}L^i_{\underset\vee{\alpha n}}
      +L^\alpha_{\underset\vee{jn}}L^i_{\underset\vee{\alpha m}}
      +2L^i_{\underset\vee{j\alpha}}L^\alpha_{\underset\vee{mn}},\label{eq:rho8}
      \\\displaybreak[0]
      &\overset1{\underset9\rho}{}^i_{jmn}=R^i_{jmn}
      +L^i_{\underset\vee{jm}|n}
      -L^i_{\underset\vee{jn}|m}
      +L^\alpha_{\underset\vee{jm}}L^i_{\underset\vee{\alpha n}}
      -L^\alpha_{\underset\vee{jn}}L^i_{\underset\vee{\alpha m}}
      +2L^i_{\underset\vee{j\alpha}}L^\alpha_{\underset\vee{mn}},\label{eq:rho9}
      \\\displaybreak[0]
      &\overset1{\underset{10}\rho}{}^i_{jmn}=R^i_{jmn}
      -L^i_{\underset\vee{jm}|n}
      +L^i_{\underset\vee{jn}|m}
      -L^\alpha_{\underset\vee{jm}}L^i_{\underset\vee{\alpha n}}
      +L^\alpha_{\underset\vee{jn}}L^i_{\underset\vee{\alpha m}}
      +2L^i_{\underset\vee{j\alpha}}L^\alpha_{\underset\vee{mn}},\label{eq:rho10}
      \\\displaybreak[0]
      &\overset1{\underset{11}\rho}{}^i_{jmn}=R^i_{jmn}
      -L^i_{\underset\vee{jm}|n}
      +L^i_{\underset\vee{jn}|m}
      +L^\alpha_{\underset\vee{jm}}L^i_{\underset\vee{\alpha n}}
      -L^\alpha_{\underset\vee{jn}}L^i_{\underset\vee{\alpha m}}
      -2L^i_{\underset\vee{j\alpha}}L^\alpha_{\underset\vee{mn}},\label{eq:rho11}
      \\\displaybreak[0]
      &\overset1{\underset{12}\rho}{}^i_{jmn}=R^i_{jmn}
      -L^i_{\underset\vee{jm}|n}
      +L^i_{\underset\vee{jn}|m}
      -L^\alpha_{\underset\vee{jm}}L^i_{\underset\vee{\alpha n}}
      +L^\alpha_{\underset\vee{jn}}L^i_{\underset\vee{\alpha m}}
      -2L^i_{\underset\vee{j\alpha}}L^\alpha_{\underset\vee{mn}},\label{eq:rho12}
      \\\displaybreak[0]
      &\overset1{\underset{13}\rho}{}^i_{jmn}=R^i_{jmn}
      +L^i_{\underset\vee{jm}|n}
      -L^i_{\underset\vee{jn}|m}
      +L^\alpha_{\underset\vee{jm}}L^i_{\underset\vee{\alpha n}}
      +L^\alpha_{\underset\vee{jn}}L^i_{\underset\vee{\alpha m}}
      ,\label{eq:rho13}
      \\\displaybreak[0]
      &\overset1{\underset{14}\rho}{}^i_{jmn}=R^i_{jmn}
      -L^i_{\underset\vee{jm}|n}
      -L^i_{\underset\vee{jn}|m}
      +L^\alpha_{\underset\vee{jm}}L^i_{\underset\vee{\alpha n}}
      +L^\alpha_{\underset\vee{jn}}L^i_{\underset\vee{\alpha m}}
      .\label{eq:rho14}
    \end{align}


\begin{thebibliography}{33}




\bibitem{blau} \textbf{M. Blau}, \emph{Lecture Notes on General
    Relativity}, arXiv: 9712019.

\bibitem{cvetkoviczlatanovic1} \textbf{M. D. Cvetkovi\'c, M. Lj. Zlatanovi\'c},
\emph{New Cartan's Tensors and Pseudotensors in a Generalized
Finsler Space}, Filomat 28:1 (2014), 107--117.

    \bibitem{e1}
  {\bf A. Einstein}, \emph{A generalization of the
  relativistic theory of gravitation}, Ann. of. Math., 45 (1945),
  No. 2, 576--584.

  \bibitem{e2}
  {\bf A. Einstein}, \emph{Bianchi identities in the
  generalized theory of gravitation}, Can. J. Math., (1950), No. 2,
  120--128.

  \bibitem{e3}
  {\bf A. Einstein}, \emph{Relativistic Theory of
  the Non-symmetric Field}, Princeton University Press, New Jersey,
  1954, 5th edition.

      \bibitem{eisNRG} \textbf{L. P. Eisenhart}, \emph{Non-Riemannian
    Geometry}, New York, 1927.

    \bibitem{eis01}
  {\bf L. P. Eisenhart}, \emph{Generalized Riemannian
  spaces}, Proc. Natl. Acad. Sci. USA 37 (1951) 311--315.

  \bibitem{eis02}
  {\bf L. P. Eisenhart}, \emph{Generalized Riemannian
  spaces}, II, Proc. Natl. Acad. Sci. USA 38 (1952) 505--508.

  \bibitem{greiner} \textbf{W. Greiner, J. Reinhardt}, \emph{Field
  Quantization}, Springer-Verlag Berlin Heidelberg New York, Berlin
  1996.



      \bibitem{z2} \textbf{S. Ivanov, M. Lj. Zlatanovi\'c},
    \emph{Conenctions on a non-symmetric (generalized) Riemannian
    manifold and gravity}, Class. Quantum Grav. 33 (2016) 075016.




        \bibitem{mik5} \textbf{J. Mike\v s, E. Stepanova,
    A. Van\v zurova, et al.}, \emph{Differential geometry of special mappings},
    Olomouc: Palacky University, 2015.


    \bibitem{mincic4} \textbf{S. M. Min\v ci\'c}, \emph{Curvature tensors of the space of non-symmetric affine connexion,
    obtained from the curvature pseudotensors}, Matemati\v cki Vesnik,
    13 (28) (1976), 421--435.

    \bibitem{mincic2} \textbf{S. M. Min\v ci\'c}, \emph{Independent
    curvature tensors and pseudotensors of spaces with
    non-symmetric affine connexion}, Coll. Math. Soc. J\'anos Bolayai,
    31. Dif. geom., Budapest (Hungary), (1979), 445--460.

    \bibitem{mincicnovi} \textbf{S. M. Min\v ci\'c}, \emph{On Ricci
    Type Identities in Manifolds With Non-Symmetric Affine
    Connection}, Publications De L'Institut Math\'ematique, Nouvelle
    s\'erie, tome 94 (108) (2013), 205--217.

    \bibitem{mincvel} \textbf{S. M. Min\v ci\'c, Lj. S.
    Velimirovi\'c}, \emph{Spaces With Non-Symmetric Affine
    Connection}, Novi Sad J. Math., Vol. 38, No. 3, 2008, 157--164.



          \bibitem{mileva1} \textbf{M. Prvanovi\'c}, \emph{Four curvature tensors of non-symmetric affine
    connexion}, (in Russian),
    Proceedings of the conference "150 years of Lobachevski geometry", Kazan' 1976,
    Moscow 1997, 199--205.

        \bibitem{shapiro} \textbf{I. L. Shapiro},
    \emph{Physical Aspects of the Space-Time Torsion},
    Physics Reports, Volume 357, Issue 2, January 2002,
    113--213.

            \bibitem{sinjukov} \textbf{N. S. Sinyukov}, \emph{Geodesic mappings
    of Riemannian spaces},  (in Russian), "Nauka", Moscow, 1979.





          \bibitem{z4} \textbf{M. Lj. Zlatanovi\'c}, \emph{New projective
    tensors for equitorsion geodesic mappings}, Applied Mathematics
    Letters 25 (2012), No. 5, 890--897.



  \end{thebibliography}
\end{document}